\documentclass[11pt,a4paper]{article}

\usepackage[a4paper,left=1.65cm,right=1.65cm,top=1.65cm,bottom=1.65cm]{geometry}
\usepackage{t1enc}
\usepackage[utf8]{inputenc}
\usepackage{amsthm,amssymb}

\usepackage[mathcal,mathscr]{eucal}
\usepackage{enumerate}
\usepackage[sort,nocompress]{cite}
\usepackage{indentfirst}
\usepackage{todonotes}
\usepackage{subcaption}
\usepackage{hyperref}
\hypersetup{
    colorlinks=true,       
    linkcolor=blue,        
    citecolor=red,         
    filecolor=magenta,     
    urlcolor=cyan,         
    linktocpage=true
}

\theoremstyle{plain}
\newtheorem{theorem}{Theorem}
\newtheorem{lemma}[theorem]{Lemma}
\newtheorem{corollary}[theorem]{Corollary}
\newtheorem{claim}[theorem]{Claim}

\theoremstyle{definition}

\newtheorem{remark}[theorem]{Remark}

\usepackage{amsmath}

\begin{document}

\title{Supermodularity in Unweighted Graph Optimization I:  \break
Branchings and Matchings }

\author{Krist\'of B\'erczi\thanks{ MTA-ELTE Egerv\'ary Research Group,
Department of Operations Research, E\"otv\"os University, P\'azm\'any
P. s. 1/c, Budapest, Hungary, H-1117. e-mail:  {\tt berkri\char'100
cs.elte.hu .}} \ \ and  \ {Andr\'as Frank\thanks{MTA-ELTE Egerv\'ary
Research
Group, Department of Operations Research, E\"otv\"os University,
P\'azm\'any P. s. 1/c, Budapest, Hungary, H-1117. e-mail:  {\tt
frank\char'100 cs.elte.hu .} } } }

\maketitle

\begin{abstract} The main result of the paper is motivated by the
following two, apparently unrelated graph optimization problems:  (A)
\ as an extension of Edmonds' disjoint branchings theorem,
characterize digraphs comprising $k$ disjoint branchings $B_i$ each
having a specified number $\mu _i$ of arcs, (B) \ as an extension of
Ryser's maximum term rank formula, determine the largest possible
matching number of simple bipartite graphs complying with
degree-constraints.  The solutions to these problems and to their
generalizations will be obtained from a new min-max theorem on
covering a supermodular function by a simple degree-constrained
bipartite graph.  A specific feature of the result is that its minimum
cost extension is already {\bf NP}-hard.  Therefore classic
polyhedral tools themselves definitely cannot be sufficient for
solving the problem, even though they make some good service in our
approach.  \end{abstract}

\section{Introduction}

Network flow theory provides a basic tool to treat conveniently
various graph characterization and optimization problems such as the
degree-constrained subgraph problem in a bipartite graph (or bigraph,
for short) or the $k$ edge-disjoint $st$-paths problem in a directed
graph (or digraph, for short).  (Throughout the paper, we use the
terms bipartite graph and bigraph as synonyms, and simlarly the terms
directed graph and digraph.)  Another general framework in graph
optimization is matroid theory.  For example, the problem of extending
$k$ given subtrees of a graph to $k$ disjoint spanning trees can be
solved with the help of matroids, as well as the problem of finding a
cheapest rooted $k$-edge-connected subgraph of a digraph.

A common generalization of these two big branches of combinatorial
optimization is the theory of submodular flows, initiated by Edmonds
and Giles \cite{Edmonds-Giles}.  This covers not only the basic
results on maximum flows and min-cost circulations from network flow
theory and weighted (poly)matroid intersection or matroid partition
from matroid theory but also helps solving significantly more complex
graph optimization problems such as the one of finding a minimum
dijoin in a digraph (the classic theorem of Lucchesi and Younger) or
finding a $k$-edge-connected orientation of a mixed graph.

However general is the framework of submodular flows, it leaves open
one of the most significant unsolved questions of matroid optimization
concerning the existence of $k$ (or just 2) disjoint common bases of
two matroids.  This is settled only in special cases, among them the
most important one is a theorem of Edmonds \cite{Edmonds73} on the
existence of $k$ disjoint spanning arborescences of common root in a
digraph.  This version is sometimes called the weak form of Edmonds'
theorem while its strong form characterizes digraphs admitting $k$
disjoint spanning branchings with prescribed root-sets.  Due to the
specific position of Edmonds' theorem within combinatorial
optimization, it is particularly important to investigate its
extensions and variations.  For example, the problem of finding $k$
disjoint spanning arborescences with no requirements on the location
of their roots is a nicely tractable version \cite{FrankP4}, and even
more generally, one may impose upper and lower bounds for each node
$v$ to constrain the number of arborescences rooted at $v$.  By using
analogous techniques, one can characterize digraphs comprising $k$
disjoint spanning branchings each having $\mu $ arcs.

A characteristic feature of submodular flows is that the corresponding
linear system is totally dual integral and therefore the weighted (or
minimum cost) versions of the graph theoretic applications are
typically also tractable.  For example, not only the minimum
cardinality dijoin problem can be solved in polynomial time but its
minimum cost version as well \cite{FrankJ7}.  Or, via submodular
flows, there is a polynomial time algorithm to find a cheapest
$k$-edge-connected orientation of a $2k$-edge-connected graph.

More generally, a great majority of min-max theorems and good
characterizations in combinatorial optimization has a polyhedral
background that makes possible to manage weighted or min-cost versions
(see, for example non-bipartite matchings).  In this view, it is quite
interesting that around the same time when submodular flows were
introduced, pretty natural graph optimization problems emerged in
which the minimum cardinality case was shown to be polynomially
solvable while the weighted version turned out to be {\bf
NP}-complete.  For example, Eswaran and Tarjan \cite{ET} found a
min-max formula and an algorithm to make a digraph strongly connected
by adding a minimum number of new arcs but the minimum cost version of
the problem is clearly {\bf NP}-complete as the directed Hamiltonian
circuit problem is a special case.  Therefore no polyhedral approach
can exist for this augmentation problem.  (Note that the original
cardinality version of Eswaran and Tarjan has nothing to do with the
problem of packing common bases of two matroids.)

\medskip

\begin{remark} The problem of deciding whether a digraph includes a
circuit of cost at most a specified number is {\bf NP}-complete.  The
closely related minimization version of this problem consists of
determining the minimum cost of a directed circuit.  Strictly
speaking, this minimization form is {\bf NP}-hard and not {\bf
NP}-complete since it is not known to be in {\bf NP}.  But the two
problems are informally so close to each other that we take the
liberty to use throughout the inaccurate term {\bf NP}-complete for
the minimizing form, as well.\end{remark}

\medskip

Recently, it turned out that the roots of a somewhat similar
phenomenon go back to as early as 1958 when Ryser \cite{Ryser58}
solved the maximum term rank problem (which is equivalent to finding a
simple bipartite graph $G$ with a specified degree sequence so that
$G$ has a matching with cardinality at least a specified number $\ell
$, or equivalently, the matching number $\nu (G)$ of $G$ is as large
as possible).  The minimum cost version of this problem had not been
settled for a long time.  Ford and Fulkerson, for example, considered
a natural attempt by using network flows but they concluded in their
book \cite{Ford-Fulkerson} that the flow approach did not seem to work
in this case.  (For the exact citation, see Section \ref{termrank}.)
Recently, however, it was shown (\cite{KaMk}, \cite{Palvolgyi},
\cite{Puleo}) that this min-cost version of the maximum term rank
problem is {\bf NP}-complete.

Therefore the failure of using network flows to attack the maximum
term rank problem was not by chance at all, and the same {\bf
NP}-completeness result shows that even submodular flows could not be
able to help.  The sharp borderline between the problem of finding a
degree-specified simple bipartite graph and the problem of finding a
degree-specified simple bipartite graph with matching number at least
$\ell$ is best clarified by the fact that --though both problems are
in {\bf P}-- the natural extension of the first problem, when a
degree-specified subgraph of an initial bipartite graph is to be
found, is still in {\bf P}, while the analogous extension of the
second problem, when a degree-specified subgraph with matching number
at least $\ell$ of an initial bipartite graph is to be found, is
already {\bf NP}-complete.

In a paper by the second author \cite{FrankJ23}, a min-max theorem was
developed to solve the general edge-connectivity augmentation problem
of digraphs.  It was shown in \cite{FrankS5}, that edge-connectivity
augmentation problem of digraphs could be embedded in an abstract
framework concerning optimal arc-covering of supermodular functions.
That min-max theorem seems to be the very first appearance of a
min-max result on sub- or supermodular functions in which the weighted
version included {\bf NP}-complete problems.

Frank and Jord\'an \cite{FrankJ31} generalized this result further and
proved a min-max theorem on optimally covering a so-called
supermodular bi-set function by digraphs.  We shall refer to the main
result of \cite{FrankJ31} (and its equivalent reformulation, too) as
the supermodular arc-covering theorem (Theorem \ref{Frank-Jordan}
below).  It should be emphasized that this framework
characteristically differs from previous models using sub- or
supermodular functions, such as polymatroids or submodular flows,
since it solves such cardinality optimization problems for which the
corresponding weighted versions are {\bf NP}-complete.  One of the
most important applications was a solution to the minimum directed
node-connectivity augmentation problem but several other problems
could be treated in this way.  For example, with its help, the
degree-sequences of $k$-edge-connected and $k$-node-connected digraphs
could be characterized (without requiring simplicity of the realizing
digraph).  Also, it implied (an extension of) Gy{\H o}ri's
\cite{Gyori85} beautiful theorem on covering a vertically convex
polyomino by a minimum number of rectangles.  Yet another application
described a min-max formula for $K_{t,t}$-free $t$-matchings of a
bipartite graph \cite{FrankJ46}.  In a recent application, Soto and
Telha \cite{Soto-Telha14} described an elegant extension of Gy{\H
o}ri's theorem.

One may consider analogous problems concerning simple digraphs
covering supermodular functions.  Unfortunately, it turned out
recently that the problem of supermodular coverings with {\it simple}
digraph includes {\bf NP}-complete special cases.  Therefore there is
no hope to develop a general version of the min-max theorem of Frank
and Jord\'an where the covering digraph is requested to be simple.

The present work is the first member of a series of four papers.  Our
general goal is to describe special cases where simplicity can
successfully be treated.  Here a new min-max theorem is developed on
covering an intersecting supermodular function with a simple
degree-constrained bipartite graph.  One application is a new theorem
on disjoint branchings which provides a necessary and sufficient
condition for the existence of $k$ disjoint spanning branchings
$B_1,\dots ,B_k$ in a digraph such that, for each $i=1,\dots ,k$, the
cardinality of $\vert B_i\vert $ lies between prescribed lower and
upper bounds $f_i$ and $g_i$, and such that, for each node $v\in V$,
the (total) in-degree $\varrho _F(v)$ of $v$ lies between specified
lower and upper bounds $f_{\rm in}(v)$ and $g_{\rm in}(v)$, where
$F=B_1\cup \cdots \cup B_k$.

As another consequence, we shall show that Ryser's maximum term rank
problem nicely fits this new framework.  Ryser's original maximum term
rank theorem, in equivalent terms of bipatite graphs, provides a
min-max formula for the maximum of the matching number of simple
bipartite graphs meeting a degree prescription.  We not only show how
this result follows from our general supermodular framework but a
significantly more difficult extension will also be derived in which
the goal is to determine the maximum of the matching number of
degree-{\rm constrained} simple bigraphs.

In Part II \cite{Berczi-Frank16b} of the series, matroidal
generalization of the new framework is described which gives rise to a
matroidal extension of Ryser's maximum term rank theorem.  We also
develop the more general augmentation version of Ryser's max term rank
formula, when some edges of the graph (correspondingly, some 1's of
the matrix) are specified.

In Part III \cite{Berczi-Frank16c}, yet another special case of the
supermodular arc-covering theorem is analysed where simplicity of the
covering digraph is tractable, and we derive there, among others, a
characterization of degree-sequences of simple $k$-node-connected
digraphs, providing in this way a straight generalization of a recent
result of Hong, Liu, and Lai \cite{Hong-Liu-Lai} on the
characterization of degree-sequences of simple strongly connected
digraphs.  The approach also gives rise to a characterization for the
augmentation problem where an initial digraph is to be augmented to
obtain a $k$-node-connected, simple, degree-specified digraph.

Part IV \cite{Berczi-Frank17} will be devoted to explore algorithmic
aspects of the problems.  Based on the ellipsoid method, there is a
polynomial time algorithm \cite{FrankJ31} to compute the optima in the
supermodular arc-covering theorem (Theorem \ref{Frank-Jordan} below).
Therefore our approach ensures polynomial algorithms for several
applications we discuss in the three papers (but not for all).  In
addition, V\'egh and Bencz\'ur \cite{Benczur-Vegh} developed a purely
combinatorial algorithm for the directed node-connectivity
augmentation problem, and their algorithm can be extended to the
general supermodular arc-covering theorem, as well, provided a
submodular function minimizing oracle is available.  This version is
polynomial in the size of the ground-set and in the maximum value of
the supermodular function to be covered (in other words, the algorithm
is pseudo--polynomial).  The algorithm of V\'egh and Bencz\'ur,
however, is pretty intricate and one goal of Part IV is to develop
simpler algorithms.  For example, a purely graph theoretical algorithm
(not relying on a submodular function minimizing oracle) will be
constructed to compute $k$ disjoint branchings $B_1,\dots ,B_k$ in a
digraph with sizes $\mu _1,\dots ,\mu _k$.  Another goal of Part IV will be
to develop an algorithmic solution to the degree-constrained
augmentation version of the maximum term rank problem, a problem
discussed in Part II.

\subsection{Notions and notation}

We close this introductory section by mentioning notions and notation.

For a number $x$, let $x\sp +:=\max\{x,0\}$.  For a function
$m:V\rightarrow {\bf R},$ the set-function $\widetilde m$ is defined
by $\widetilde m(X)=\sum [m(v):v\in X]$.  A set-function $p$ can
analogously be extended to families $\cal F$ of sets by $\widetilde
p({\cal F})= \sum [p(X):X\in {\cal F}]$.

Two subsets $X$ and $Y$ of a ground-set $V$ are {\bf comparable} if
$X\subseteq Y$ or $Y\subseteq X$, {\bf intersecting} if $X\cap Y\not
=\emptyset $, {\bf properly intersecting} if they are non-comparable
and intersecting, {\bf crossing} if none of the sets $X-Y, Y-X, X\cap
Y, V-(X\cup Y)$ is empty.

For two non-empty subsets $S$ and $T$ of $V$, the subsets $X,Y$ are
{\bf $ST$-independent} if $X\cap Y\cap T=\emptyset $ or $S-(X\cup
Y)=\emptyset $, {\bf $ST$-crossing} if they are non-comparable, $X\cap
Y\cap T\not =\emptyset $, and $S-(X\cup Y)\not =\emptyset $. \ $X$ and
$Y$ are {\bf $T$-intersecting} if $X\cap Y\cap T\not =\emptyset $, and
{\bf properly $T$-intersecting} if they are non-comparable and $X\cap
Y\cap T\not =\emptyset $. Typically, we do not distinguish between a
one-element set $\{v\}$, called a {\bf singleton}, and its only
element $v$.

For an arc $f=uv$, node $v$ is the {\bf head} of $f$ and $u$ is its
{\bf tail}.  The arc $uv$ {\bf enters} or {\bf covers} a subset
$X\subset V$ if $u\in V-X$ and $v\in X$.  Given a digraph $D=(V,A)$,
the {\bf in-degree} of a subset $X\subseteq V$ is the number of arcs
entering $X$, denoted by $\varrho _D(X)$ or $\varrho _A(X)$.  The {\bf
out-degree} $\delta _D(X)=\delta _A(X)$ is the number of arcs leaving
$X$.  An arc $st$ is an {\bf $ST$-arc} if $s\in S$ and $t\in T$.  A
digraph $D=(V,A)$ {\bf covers} a set-function $p$ on $V$ if $\varrho
_D(X)\geq p(X)$ holds for every subset $X\subseteq V$.

An arc with coinciding head and tail is called a {\bf loop}.  Two arcs
from $s$ to $t$ are called {\bf parallel}.  A digraph with no loops
and parallel arcs is {\bf simple}.  Note, however, that simple
digraphs are allowed to have two oppositely oriented arcs $uv$ and
$vu$.  Simplicity of an undirected graph is defined analogously.

Let $G=(S,T;E)$ be a bipartite graph.  For a subset $Y\subseteq T$,
let $$ \hbox{ $\Gamma _G(Y)=\{s\in S:$ there is an edge $st\in E$ with
$t\in Y\}$, }\ $$ that is, $\Gamma _G(Y)$ is the set of neighbours of
$Y$.  We say that $G$ {\bf covers} a set-function $p_T$ on $T$ if \begin{equation}
\vert \Gamma _G(Y)\vert \geq p_T(Y) \ \hbox{for every subset
$Y\subseteq T$.  }\ \label{(bigraph.cover)} \end{equation}

Even if it is not mentioned explicitly, we assume throughout that each
set-function is zero on the empty set.  Also, the empty sum is defined
to be zero.  A set-function $p$ on $T$ is {\bf monotone
non-decreasing} if $p(X)\leq p(Y)$ whenever $\emptyset \subset
X\subseteq Y\subseteq T$.  For a set-function $b$ on ground-set $V$,
\begin{equation}b(X)+b(Y)\geq b(X\cap Y) + b(X\cup Y) \label{(submod)} \end{equation} is
called the {\bf submodular inequality} on $X,Y\subseteq V$.

The function $b$ is {\bf fully} (respectively, {\bf intersecting,
crossing}) {\bf submodular} if \eqref{(submod)} holds for each (resp.,
intersecting, crossing) sets $X$ and $Y$.  Fully submodular functions
will often be mentioned simply as submodular.  A set-function $p$ is
{\bf supermodular} if $-p$ is submodular, {\bf positively
intersecting} ({\bf crossing, $ST$-crossing}) supermodular if the
supermodular inequality $$p(X)+p(Y)\leq p(X\cap Y)+p(X\cup Y)$$ holds
for intersecting (crossing, $ST$-crossing) subsets for which $p(X)>0$
and $p(Y)>0.$ The {\bf complementary function} $p$ of a set-function
$b$ with finite $b(V)$ is defined by $$p(X):=b(V)-b(V-X).$$ Clearly,
$b$ is submodular if and only if $p$ is supermodular.  For a pair
$(p,b)$ of set-functions, \begin{equation}b(X)-p(Y)\geq b(X-Y) - p(X\cup Y)
\label{(supermod)} \end{equation} is called the {\bf cross-inequality} on
$X,Y\subseteq V$.  The pair is called {\bf paramodular} ({\bf
intersecting paramodular}) if $b$ is (intersecting) submodular, $p$ is
(intersecting) supermodular and the cross-inequality holds for every
(properly intersecting) $X$ and $Y$.  For a paramodular pair $(p,b)$,
the polyhedron $$Q(p,b)=\{x\in {\bf R}\sp V:  p\leq \widetilde x\leq
b\}$$ is called a {\bf generalized polymatroid} or {\bf
g-polymatroid}.  By convention, the empty set is also considered to be
a g-polymatroid.  For a submodular function $b$ with $b(V)$ finite,
the polyhedron $B(b):=\{x\in {\bf R}\sp V:  \ \widetilde x\leq b, \
\widetilde x(V)=b(V)\}$ is called a {\bf base-polyhedron} and we speak
of a {\bf 0-base-polyhedron} if $b(V)=0$.  For a supermodular function
$p$ with finite $p(V)$, the polyhedron $B'(p):=\{x\in {\bf R}\sp V:  \
\widetilde x\geq p, \ \widetilde x(V)=p(V)\}$ is also a
base-polyhedron since $B'(p)=B(b)$ holds for the complementary
function $b$ of $p$.

\medskip

All the notions, notation, and terminology not mentioned explicitly in
the paper can be found in the book of the second author
\cite{Frank-book}.

\medskip \medskip

\section{Background results}

\subsection{Degree-specified and degree-constrained bipartite graphs}

\subsubsection{Subgraph problems}

Let $S$ and $T$ be two disjoint sets and $V:=S\cup T$.  Our starting
point is the classic Hall theorem:

\begin{theorem}\label{Hall} A bigraph $G=(S,T;E)$ has a matching covering
$T$ if and only if \begin{equation}\vert \Gamma _G(Y)\vert \geq \vert Y\vert \
\hbox{for every subset $Y\subseteq T$.  }\ \label{(Hall.felt)} \end{equation}
$G$ has a perfect matching if and only if $\vert S\vert =\vert T\vert
$ and \eqref{(Hall.felt)} holds.  \end{theorem}

For a given non-negative integer-valued function $m:V\rightarrow {\bf
Z}_+$, its restrictions to $S$ and to $T$ are denoted by $m_S$ and
$m_T$, respectively.  We also use the notation $m=(m_S,m_T)$.  It is
assumed throughout that $\widetilde m_S(S)=\widetilde m_T(T)$ and this
common value will be denoted by $\gamma $. We say that $m$ or the pair
$(m_S,m_T)$ is a {\bf degree-specification} and that a bipartite graph
$G=(S,T;E)$ {\bf fits} or {\bf meets} this degree-specification if
$d_G(v)=m(v)$ holds for every node $v\in V$.

\begin{theorem}[Ore \cite{Ore56}] \label{basic.GRalt} Let $G_0=(S,T;E_0)$ be
a bipartite graph and $m=(m_S,m_T)$ a degree-specification for which
$\widetilde m_S(S)=\widetilde m_T(T)=\gamma $. There is a subgraph
$G=(S,T;E)$ of $G_0$ fitting the degree-specification $m$ if and only
if \begin{equation}\widetilde m_S(X) + \widetilde m_T(Y) -d_{G_0}(X,Y) \leq \gamma
\ \hbox{whenever $X\subseteq S, \ Y\subseteq T$}\
\label{(basic.GRalt)} \end{equation} where $d_{G_0}(X,Y)$ denotes the number of
edges connecting $X$ and $Y$.  \end{theorem}

Let $g_S:S\rightarrow {\bf Z}_+$ and $g_T:T\rightarrow {\bf Z}_+$ be
upper bound functions while $f_S:S\rightarrow {\bf Z}_+$ and
$f_T:T\rightarrow {\bf Z}_+$ lower bound functions.  Let
$f_V=(f_S,f_T)$ and $g_V=(g_S,g_T)$ and assume that $f_V\leq g_V$.
Call a bipartite graph $G=(S,T;E)$ {\bf $(f_T,g_S)$-feasible} if \begin{equation}
\hbox{ $d_G(s)\leq g_S(s)$ for every $s\in S$ and $d_G(t)\geq f_T(t)$
for every $t\in T.$ }\ \label{(ftgs0)} \end{equation} \noindent The bigraph
$G=(S,T;E)$ (and its degree function $d_G$) is said to be {\bf
$(f_V,g_V)$-feasible} or {\bf degree-constrained by} $(f_V,g_V)$ if
$f_V(v)\leq d_G(v)\leq g_V(v)$ holds for every node $v\in V$.

\begin{theorem}[Linking property, Ford and Fulkerson] \label{basic.linking}
Let $G_0=(S,T;E_0)$ be a bipartite graph.  Let $g_S:S\rightarrow {\bf
Z}_+$ and $g_T:T\rightarrow {\bf Z}_+$ be upper bound functions while
$f_S:S\rightarrow {\bf Z}_+$ and $f_T:T\rightarrow {\bf Z}_+$ lower
bound functions.  There is an $(f_V,g_V)$-feasible subgraph $G$ of
$G_0$ if and only if there is an $(f_S,g_T)$-feasible subgraph $G'$ of
$G_0$ and there is an $(f_T,g_S)$-feasible subgraph $G''$ of $G_0$.
\end{theorem}

With standard techniques, such as network flows or total
unimodularity, the following theorem can also be derived.

\begin{theorem}\label{ujintro.link2} Suppose that a bigraph $G_0$ has a
subgraph degree-constrained by $(f_V,g_V)$.  \ $G_0$ has an
$(f_V,g_V)$-feasible subgraph $G=(S,T;E)$:

\noindent {\bf (A)} \ for which $\alpha \leq \vert E\vert $ if and
only if \begin{equation}\hbox{ $\widetilde g_S(S-X) + \widetilde g_T(T-Y) +
d_{G_0}(X,Y) \geq \alpha $ whenever $X\subseteq S, \ Y\subseteq T$, }\
\label{(ujintro.galfa0)} \end{equation}

\noindent {\bf (B)} \ for which $\vert E\vert \leq \beta $ if and only
if \begin{equation}\hbox{ $\widetilde f_S(X) + \widetilde f_T(Y) - d_{G_0}(X,Y)
\leq \beta $ \ whenever $X\subseteq S, \ Y\subseteq T$, }\
\label{(ujintro.fbeta0)} \end{equation}

\noindent {\bf (AB)} \ for which $\alpha \leq \vert E\vert \leq \beta
$ if and only if both \eqref{(ujintro.galfa0)} and
\eqref{(ujintro.fbeta0)} hold.\end{theorem}

It should be noted that the \lq only if\rq \ part in the theorems
above are straightforward.  For example, in Theorem \ref{basic.GRalt},
we can argue that in a subgraph $G$ of $G_0$ fitting $m$ there are
$\widetilde m_S(X)$ edges leaving $X$ and $\widetilde m_T(Y)$ edges
leaving $Y$, and since at most $d_{G_0}(X,Y)$ edges may be counted
twice, the total number $\gamma $ of edges of $G$ is at least
$\widetilde m_S(X) + \widetilde m_T(Y) -d_{G_0}(X,Y)$.

\subsubsection{Synthesis problems}

When the initial graph $G_0$ is the complete bipartite graph on $S$
and $T$, the theorems can be simplified.  Let ${\cal G}(m_S,m_T)$
denote the set of simple bipartite graphs fitting $(m_S,m_T)$.  Gale
\cite{Gale57} and Ryser \cite{Ryser57} found, in an equivalent form,
the following characterization.

\begin{theorem}[Gale and Ryser] \label{ujterm.GR} There is a simple
bipartite graph $G$ fitting the degree-specification $m$ if and only
if \begin{equation}\widetilde m_S(X) + \widetilde m_T(Y) -\vert X\vert \vert
Y\vert \leq \gamma \ \hbox{whenever $X\subseteq S, \ Y\subseteq T.$}\
\label{(ujterm.GR)} \end{equation} Moreover, \eqref{(ujterm.GR)} holds if the
inequality is required only when $X$ consists of the $i$ elements of
$S$ having the $i$ largest values of $m_S$ and $Y$ consists of the $j$
elements of $T$ having the $j$ largest values of $m_T$ $(i=1,\dots
,\vert S\vert , \ j=1,\dots ,\vert T\vert )$.  \end{theorem}

Instead of exact degree-specifications, one may impose upper and/or
lower bounds for the degrees.

\begin{theorem}\label{ujterm.gf} Let $g_S:S\rightarrow {\bf Z}_+$ be an
upper bound function on $S$ and let $f_T:T\rightarrow {\bf Z}_+$ be a
lower bound function on $T$.  There is an $(f_T,g_S)$-feasible simple
bipartite graph $G$ if and only if \begin{equation}\widetilde g_S(X) + \widetilde
f_T(Y) -\vert X\vert \vert Y\vert \leq \widetilde g_S(S) \
\hbox{whenever $X\subseteq S, \ Y\subseteq T.$ }\ \label{(ujterm.gf)}
\end{equation} Moreover, \eqref{(ujterm.gf)} holds if the inequality is
required only when $X$ consists of elements with the $i$ largest
values of $m_S$ and $Y$ consists of elements with the $j$ largest
values of $m_T$ $(i=1,\dots ,\vert S\vert , \ j=1,\dots ,\vert T\vert
)$.  \end{theorem}

The linking property formulated in Theorem \ref{basic.linking} can
also be specialized to the case when $G_0$ is the complete bipartite
graph $G\sp *=(S,T;E\sp *)$.

\begin{theorem}\label{basic.linking.spec} If there is a simple
$(f_T,g_S)$-feasible bipartite graph and there is a simple
$(f_S,g_T)$-feasible bipartite graph, then there is a simple
$(f_V,g_V)$-feasible bipartite graph.  \end{theorem}

When $G_0$ is the complete bigraph on $S\cup T$, Theorem
\ref{ujintro.link2} specializes to the following synthesis-type
problem.

\begin{theorem}\label{ujintro.link2.sint} Suppose that there is a simple
bigraph degree-constrained by $(f_V,g_V)$.  There is a simple bigraph
degree-constrained by $(f_V,g_V)$:

\noindent {\bf (A)} \ for which $\alpha \leq \vert E\vert $ if and
only if \begin{equation}\hbox{ $\widetilde g_S(S-X) + \widetilde g_T(T-Y) + \vert
X\vert \vert Y\vert \geq \alpha $ whenever $X\subseteq S, \ Y\subseteq
T,$ }\ \label{(ujintro.galfa)} \end{equation}

\noindent {\bf (B)} \ for which $\vert E\vert \leq \beta $ if and only
if \begin{equation}\hbox{ $\widetilde f_S(X) + \widetilde f_T(Y) - \vert X\vert
\vert Y\vert \leq \beta $ \ whenever $X\subseteq S, \ Y\subseteq T$,
}\ \label{(ujintro.fbeta)} \end{equation}

\noindent {\bf (AB)} \ for which $\alpha \leq \vert E\vert \leq \beta
$ if and only if both \eqref{(ujintro.galfa)} and
\eqref{(ujintro.fbeta)} hold.\end{theorem}

\subsubsection{Synthesis versus subgraph problems}

The synthesis problem of degree-constrained and degree-specified
simple bigraphs is just a special case of the corresponding subgraph
problems.  It turns out, however, that several other synthesis
problems cannot be attacked in this way since the more general
subgraph problem is already {\bf NP}-complete.  For example, it is
trivial to decide if there is a bigraph $G=(S,T;E)$ which is connected
and meets the identically 2 degree-specification since such a graph is
just a bipartite Hamiltonian circuit, and therefore the only
requirement is $\vert S\vert =\vert T\vert \geq 2$.  On the other
hand, it is known to be {\bf NP}-complete to decide if an initial
bigraph $G_0$ includes a Hamiltonian circuit.

At other occasions the situation is more complicated.  For example,
one may consider the synthesis problem of finding a simple, perfectly
matchable, and degree-specified bigraph.  This problem is solvable but
its subgraph version where a perfectly matchable degree-specified
subgraph of an initial bigraph $G_0$ has to be found is already {\bf
NP}-complete (\cite{KaMk}, \cite{Palvolgyi}, \cite{Puleo}).

\subsection{Covering supermodular functions with digraphs and
bigraphs}

\subsubsection{Covering by bigraphs}

We call a set-function $p$ on a ground-set $T$ {\bf
element-subadditive} if $p(Y)+p(t)\geq p(Y+t)$ holds whenever
$Y\subseteq T$ and $t\in T$.  The following early result on bipartite
graphs and supermodular functions is due to Lov\'asz \cite{Lovasz70a}.

\begin{theorem}\label{intro.Lovasz} Let $G_0=(S,T;E_0)$ be a simple
bipartite graph and $p_T$ a positively intersecting supermodular
function on $T$ which is, in addition, element-subadditive.  There is
a subgraph $G$ of $G_0$ covering $p_T$ for which $d_G(t) = p_T(t)$
whenever $t\in T$ if and only if \begin{equation}\vert \Gamma _{G_0}(Y)\vert \geq
p_T(Y) \ \hbox{holds for every subset}\ Y\subseteq T.
\label{(intro.Lovasz.felt)} \end{equation} \end{theorem}

This was extended by Frank and Tardos \cite{FrankJ19} as follows.

\begin{theorem}\label{Frank-Tardos} Let $G_0=(S,T;E_0)$ be a simple
bipartite graph and $p_T$ a positively intersecting supermodular
function on $T$.  Let $g_T:T\rightarrow {\bf Z}_+$ be an upper bound
function.  There is a subgraph $G$ of $G_0$ covering $p_T$ for which
$d_G(t) \leq g_T(t)$ whenever $t\in T$ if and only if \begin{equation}\hbox{
$\vert \Gamma _{G_0}(Z)\vert \geq p_T(Y\cup Z) - \widetilde g_T(Y)$
holds for disjoint subsets $Y, Z \subseteq T$.  }\
\label{(felso.Fra-Tar.felt1)} \end{equation} \end{theorem}

It should be noted that the problem in Theorem \ref{intro.Lovasz} can
be formulated as a matroid intersection problem while the problem in
Theorem \ref{Frank-Tardos} can be cast into the submodular flow
framework.  Therefore the minimum cost versions of both cases are also
tractable.  However, both problems become {\bf NP}-complete if there
is an upper-bound $g_S$, as well, for the degrees of $G$ in $S$.

\subsubsection{Covering by digraphs}

Let $p$ be a positively $ST$-crossing supermodular function.  A basic
tool in our investigations is the following general result of Frank
and Jord\'an \cite{FrankJ31}.

\begin{theorem}[Supermodular arc-covering, set-function version]
\label{Frank-Jordan} A positively $ST$\- -cross\-ing supermodular
set-function $p$ can be covered by $\gamma $ $ST$-arcs if and only if
$\widetilde p({\cal I})\leq \gamma $ holds for every $ST$-independent
family $\cal I$ of subsets of $V$.  There is an algorithm, which is
polynomial in $\vert S\vert +\vert T\vert $ and the maximum value of
$p(X)$, to compute the minimum number of $ST$-arcs to cover $p$ and an
$ST$-independent family $\cal I$ of subsets maximizing $\widetilde
p({\cal I})$.  \end{theorem}

The theorem can be used \cite{FrankJ31} to describe characterizations
for the existence of degree-specified (and even degree-constrained)
digraphs covering $p$.  It has great many applications in graph
optimization and it serves as the major tool for the present work.  It
significantly differs from the framework of Lov\'asz above (or from
submodular flows) in that its min-cost version includes {\bf
NP}-complete special cases such as the directed Hamiltonian circuit
problem.

The existing applications give rise to a natural demand to develop a
variation of Theorem \ref{Frank-Jordan} in which no parallel arcs of
the covering digraph are allowed.  Unfortunately, this is hopeless
since the general problem includes {\bf NP}-complete special cases, as
we point out in the next theorem.  This fact underpins the
significance and the difficulties of the present work that explores
special cases of Theorem \ref{Frank-Jordan} where simplicity can be
involved.

\begin{theorem}\label{simple.npc} {\bf (A)} \ It is {\bf NP}-complete to
decide for two given degree specifications $m'\leq m$ on $V=S\cup T$
whether there exists a simple bigraph $G$ fitting $m$ which includes a
subgraph fitting $m'$.

{\bf (B)} \ The problem in Part (A) can be formulated as a special
case of the problem of finding a minimal simple digraph covering an
$ST$-crossing supermodular function.\end{theorem}

\proof{Proof.} {\bf (A)} \ By choosing $m''=m-m'$, Part (A) follows
immediately from the following elegant {\bf NP}-completeness result of
D\"urr, Guinez, and Matamala \cite{DGM}.

\begin{lemma}\label{DGMa} It is {\bf NP}-complete to decide whether, given
two degree-specifications $m'=(m'_S,m'_T)$ and $m''=(m''_S,m''_T)$,
there is a simple bigraph $G=(S,T;E)$ which can be partitioned into
two subgraphs $G'=(S,T;E')$ and $G''=(S,T;E'')$ so that $G'$ fits $m'$
and $G''$ fits $m''$.  \end{lemma}

{\bf (B)} \ Consider Theorem \ref{basic.GRalt} with $m'$ in place of
$m$.  This can be restated as follows.

\begin{claim}\label{subgraph2} A bipartite graph $G=(S,T;E)$ admits a
subgraph $G'$ fitting $m'$ if and only if \begin{equation}\varrho _D(X\cup Y) \geq
\widetilde m'_T(Y) - \widetilde m'_S(X) \ \hbox{whenever $X\subseteq
S, \ Y\subseteq T$}\ \label{(npteljes.1b)} \end{equation} where $D$ is the
digraph arising from $G$ by orienting each arc from $S$ toward
$T$.\end{claim}

Let $\gamma ':=\widetilde m'_S(S)= \widetilde m'_T(T)$ and $\gamma
:=\widetilde m_S(S)= \widetilde m_T(T)$, and define a set-function $p$
on $V$ as follows.  \begin{equation}p(V'):= \begin{cases}
\widetilde m'_T(Y) - \widetilde m'_S (X) & \text{if $Y\subseteq T, X\subseteq S, \ V'=X\cup Y, \ 1<\vert
V'\vert <\vert V\vert -1$}\\ m_T(t) & \text{if $V'=\{t\}$ for some $t\in
T$}\\ m_S(s) & \text{if $V'=V-s$ for some $s\in S$}\end{cases}
\end{equation}

\begin{claim}The set-function $p$ is $ST$-crossing supermodular.\end{claim}

\proof{Proof.} Consider the set-function $p'$ defined by $p'(V'):= \widetilde
m'_T(Y) - \widetilde m'_S (X)$ for \ $Y\subseteq T, \ X\subseteq S, \
V'=X\cup Y$.  Clearly, $p'$ is modular.  Furthermore,
$p'(t)=m'_T(t)\leq m_T(t)=p(t)$ and $p'(V-s) =\widetilde m'_T(T) -
\widetilde m'_S(S-s) = m'(s) \leq m_S(s)=p(V-s)$.  Therefore, $p$
arises from a modular function by lifting its values on elements $t$
of $T$ and by lifting its values on the complement $V-s$ of elements
$s$ of $S$.  This implies that $p$ is indeed $ST$-crossing
supermodular.  $\bullet$\endproof \medskip

Obviously, there is a simple digraph $D=(V,A)$ consisting of $\gamma $
$ST$-arcs covering $p$ if and only if there exists a simple bigraph
$G=(S,T;E)$ fitting $m$ so that \eqref{(npteljes.1b)} holds.  By Claim
\ref{subgraph2}, \eqref{(npteljes.1b)} in turn is equivalent to the
solvability of the problem in Part (A).  $\bullet $ $\bullet $\endproof

\medskip 

\section{Simple degree-specified bipartite graphs covering supermodular functions} \label{ujbasic}

Let $p_T$ be a set-function on $T$.  Recall that a bipartite graph
$G=(S,T;E)$ is said to cover $p_T$ if \begin{equation}\vert \Gamma _G(Y)\vert \geq
p_T(Y) \ \hbox{for every subset $Y\subseteq T$.  }\
\label{(basic.kell.free)} \end{equation} For example, if $p_T(Y)=\vert Y\vert $
$(Y\subseteq T)$, then \eqref{(basic.kell.free)} is the Hall-condition
\eqref{(Hall.felt)}.  Therefore Hall's theorem implies that $G=(S,T;E)$
covers $p_T$ if and only if $G$ has a matching covering $T$.  Another
special case is when $p_T(Y):=\vert Y\vert +1$\ $(\emptyset \subset
Y\subseteq T)$.  By a theorem of Lov\'asz \cite{Lovasz70a}, a bigraph
$G=(S,T;E)$ covers this $p_T$ if and only if $G$ has a forest in which
the degree of every node in $T$ is 2. This result is a direct
consequence of Theorem \ref{intro.Lovasz}.

We are interested in finding simple bipartite graphs covering $p_T$
which meet some degree-constraints (that is, upper and lower bounds)
or exact degree-specifications.  If no such constraints are imposed at
all, then the existence of a bigraph covering $p_T$ is obviously
equivalent to the requirement that

\begin{equation}\hbox{ $p_T(Y)\leq \vert S\vert $ for every $Y\subseteq T$.  }\
\label{(pS)} \end{equation}

Indeed, this condition is clearly necessary and it is also sufficient
as the complete bipartite graph $G\sp *=(S,T;E\sp *)$ covers a
set-function $p_T$ meeting \eqref{(pS)}.  Therefore we suppose
throughout that \eqref{(pS)} holds.

Our plan is the following.  First we characterize the situation when
there is a degree-prescription only on $S$.  This is then used to
settle the case when a degree-specification $(m_S,m_T)$ is given on
the whole node-set $V=S\cup T$.  In Section \ref{base-polyhedron},
with the help of a novel construction, we introduce a base-polyhedron
$B$ and prove that $(m_S,m_T)$ is realizable by a simple bigraph
covering $p_T$ precisely if the associated vector $(m_S,-m_T)$ is in
$B$.  As the intersection of a base-polyhedron with a box and with a
plank is a g-polymatroid whose non-emptiness is characterized in the
literature, this result can finally be used to handle upper and lower
bounds on the degrees of $G$ and on its edge-number.

\subsection{Degree-specification on \texorpdfstring{$S$}{S}}

Our first goal is to characterize the situation when there is a
degree-specification only on $S$.

\begin{theorem}\label{basic.base.free} Let $m_S$ be a degree-specification
on $S$ for which $\widetilde m_S(S)= \gamma $. Let $p_T$ be a
positively intersecting supermodular function on $T$ with
$p_T(\emptyset )=0$.  Suppose that \begin{equation}m_S(s)\leq \vert T\vert \
\hbox{for every $s\in S$.}\ \label{(basic.msb.free)} \end{equation} The
following statements are equivalent.  \medskip

\noindent {\bf (A)} \ There is a simple bipartite graph $G=(S,T;E)$
covering $p_T$ and fitting the degree-specification $m_S$.  \medskip

\noindent {\bf (B1)} \begin{equation}\widetilde m_S(X) + \widetilde p_T({\cal T})
- \vert {\cal T}\vert \vert X\vert \leq \gamma \ \hbox{for every
subset $X\subseteq S$ and subpartition ${\cal T}$ of $T$.  }\
\label{(basic.basefelt.1.free)} \end{equation}

\medskip \noindent {\bf (B2)} \begin{equation}\sum _{i=1}\sp q p_T(T_i) \leq \sum
_{s\in S} \min \{m_S(s), q\} \ \hbox{for every subpartition ${\cal
T}=\{T_1,\dots ,T_q\}$ of $T$.  }\ \label{(basic.felt.free)} \end{equation} \end{theorem}

\proof{Proof.} {\bf (A)} $\Rightarrow $ {\bf (B1)} \ Suppose that there is a
simple bipartite graph $G$ meeting \eqref{(basic.kell.free)}.  By the
simplicity of $G$, there are at most $\vert X\vert \vert Y\vert $
edges between $X\subseteq S$ and $Y\subseteq T$.  We claim that the
number $d_G(T_i,S-X)$ of edges between $T_i\in {\cal T}$ and $S-X$ is
at least $p_T(T_i)-\vert X\vert $. Indeed, \begin{equation}p_T(T_i) \leq \vert
\Gamma _G(T_i)\vert = \vert \Gamma _G(T_i)\cap X\vert + \vert \Gamma
_G(T_i)-X\vert \leq \vert X\vert + d_G(T_i,S-X), \label{(pTdG)} \end{equation}
that is, $d_G(T_i,S-X)\geq p_T(T_i) -\vert X\vert $. Therefore the
total number $\gamma $ of edges is at least $\widetilde m_S(X) + \sum
_i [p_T(T_i) -\vert X\vert ]$ from which
\eqref{(basic.basefelt.1.free)} follows.  \medskip

\noindent {\bf (B1)} $\Rightarrow $ {\bf (B2)} \ Let ${\cal
T}=\{T_1,\dots ,T_q\}$ be a subpartition of $T$ and let $X:=\{s\in S:
m_S(s) > q\}$.  Then \eqref{(basic.basefelt.1.free)} implies $$ \sum
_{s\in S} \min \{m_S(s), q\} = \sum [m_S(s):  s\in S-X] + q\vert
X\vert = \widetilde m_S(S-X) + q\vert X\vert =$$ $$ \gamma -\widetilde
m_S(X) + q\vert X\vert \geq \widetilde p_T({\cal T})= \sum _{i=1}\sp q
p_T(T_i), $$ as required in \eqref{(basic.felt.free)}.

\medskip

\noindent {\bf (B2)} $\Rightarrow $ {\bf (B1)} \ Let $X$ be a subset
of $S$ and ${\cal T}=\{T_1,\dots ,T_q\}$ a subpartition of $T$.  By
\eqref{(basic.felt.free)}, we have

$$ \sum _{i=1}\sp q p_T(T_i) \leq \sum _{s\in S} \min \{m_S(s), q\} =
\sum _{s\in X} \min \{m_S(s), q\} + \sum _{s\in S-X} \min \{m_S(s),
q\} \leq $$ $$ \sum _{s\in X} q + \sum _{s\in S-X} m_S(s) = q\vert
X\vert + \widetilde m_S(S-X) = q\vert X\vert + \gamma - \widetilde
m_S(X),$$ and \eqref{(basic.basefelt.1.free)} follows.  \medskip

\noindent {\bf (B1)} $\Rightarrow $ {\bf (A)} \ First of all, observe
that Condition \eqref{(basic.basefelt.1.free)}, when applied to $X:=S$
and ${\cal T}=:\{Y\}$, specializes to $\widetilde m_S(S) + p_T(Y) -
\vert S\vert \leq \gamma $, that is, \eqref{(basic.basefelt.1.free)}
implies Condition \eqref{(pS)}, and this is why explicit mentioning of
\eqref{(pS)} among the necessary conditions was avoidable.  The
following simple observation indicates that we need not concentrate on
the simplicity of $G$.

\begin{claim}\label{simplify} If there is a not-necessarily simple bipartite
graph $G=(S,T;E)$ covering $p_T$ for which $d_G(s)\leq \vert T\vert $
for each $s\in S$, then there is a simple bipartite graph $H$ covering
$p_T$ for which $d_G(s)=d_H(s)$ for each $s\in S$.  \end{claim}

\proof{Proof.} Suppose $G$ has two parallel edges $e$ and $e'$ connecting $s$
and $t$ for some $s\in S$ and $t\in T$.  Since $d_G(s)\leq \vert
T\vert $, there is a node $t'\in T$ which is not adjacent with $s$.
By replacing $e'$ with an edge $st'$, we obtain another bipartite
graph $G'$ for which $\Gamma _{G'}(Y)\supseteq \Gamma _{G}(Y)$ for
each $Y\subseteq T$, $d_{G'}(s)= d_G(s)$ for each $s\in S$, and the
number of parallel edges in $G'$ is smaller than in $G$.  By repeating
this procedure, finally we arrive at a requested simple graph.
$\bullet$\endproof \medskip

A subset $V'$ of $V:=S\cup T$ is {\bf $ST$-trivial} if no $ST$-arc
enters it, which is equivalent to requiring that $T\cap V'=\emptyset $
or $S\subseteq V'$.  We say that a subset $V'\subseteq V$ is {\bf fat}
if $V'=V-s$ for some $s\in S$ (that is, there are $\vert S\vert $ fat
sets).  The non-fat subsets of $V$ will be called {\bf normal}.  An
$ST$-independent family $\cal I$ of subsets is {\bf strongly
$ST$-independent} if any two of its normal members are
$T$-independent, that is, the intersections of the normal members of
$\cal I$ with $T$ form a subpartition of $T$.

Define a set-function $p_0$ on $V$ by \begin{equation}\hbox{ $p_0(V')=p_T(Y)-\vert
X\vert $ where $V'=X\cup Y$ whenever $X\subseteq S$ and $Y\subseteq
T$.}\ \label{(p0.def)} \end{equation} Note that $p_0$ is positively
$T$-intersecting supermodular since if $p_0(V')$ is positive, then so
is $p_T(Y)$.  Furthermore, as \eqref{(basic.basefelt.1.free)} was shown
above to imply \eqref{(pS)}, $p_0(V')=p_T(Y)-\vert X\vert $ can be
positive only if $X\not =S$ and $Y\not =\emptyset $, that is, when
$V'$ is not $ST$-trivial.

\begin{claim}\label{ujterm.mspr} $m_S(s) \geq p_0(V-s)$ holds for every
$s\in S$.  \end{claim}

\proof{Proof.} By applying \eqref{(basic.basefelt.1.free)} to $X=S-s$ and
${\cal T} = \{T\}$, we obtain that $m_S(s) \geq p_T(T)-\vert S-s\vert
= p_0(V-s)$.  $\bullet$\endproof

\medskip

Define a set-function $p_1$ on $V$ by modifying $p_0$ so as to lift
its value on fat subsets $V-s$ from $p_0(V-s)$ to $m_S(s)$ \ ($s\in
S$), that is, \begin{equation}p_1(V'):= \begin{cases}
m_S(s) & \text{if $V'=V-s$ for some
$s\in S$},\\ p_0(V') & \text{otherwise.}\end{cases} \label{(p1.def)} \end{equation}

Note that the supermodular inequality \begin{equation}p_1(V_1)+p_1(V_2)\leq
p_1(V_1\cap V_2) + p_1(V_1\cup V_2) \label{(p1.szuper)} \end{equation} holds for
$T$-intersecting normal sets with $p_1(V_1)>0$ and $p_1(V_2)>0$.

By Claim \ref{ujterm.mspr}, $p_1\geq p_0$.  In order to use Theorem
\ref{Frank-Jordan}, observe that, as $p_0$ is positively
$T$-intersecting supermodular, $p_1$ is positively $ST$-crossing
supermodular.  Let $\nu _1$ denote the maximum total $p_1$-value of a
family of $ST$-independent sets.  We call a family attaining the
maximum a {\bf $p_1$-optimizer.}

\begin{claim}\label{vanreduced} If $\cal I$ is a $p_1$-optimizer of minimum
cardinality, then $\cal I$ is strongly $ST$-independent.  \end{claim}

\proof{Proof.} Clearly, $p_1(V')\geq 0$ for each $V' \in \cal I$ for otherwise
$\cal I$ would not be a $p_1$-optimizer.  Moreover, $p_1(V')>0$ also
holds for if we had $p_1(V')=0$, then ${\cal I}-\{V'\}$ would also be
a $p_1$-optimizer contradicting the minimality of ${\cal I}$.

Suppose indirectly that $\cal I$ has two properly $T$-intersecting
normal members $V_1$ and $V_2$.  Then \eqref{(p1.szuper)} holds and
$V_1\cap V_2$ is obviously normal.  Since $\cal I$ is
$ST$-independent, we must have $S\subseteq V_1\cup V_2$ implying that
$V_1\cup V_2$ is also normal.  The inclusion $S\subseteq V_1\cup V_2$
also shows that $$p_1(V_1\cup V_2)= p_0(V_1\cup V_2) = p_T(T\cap
(V_1\cup V_2))-\vert S\vert \leq 0,$$ where the last inequality
follows from \eqref{(pS)} (which was shown above to be a consequence of
\eqref{(basic.basefelt.1.free)}).  Hence $$p_1(V_1)+p_1(V_2)\leq
p_1(V_1\cap V_2)+ p_1(V_1\cup V_2) \leq p_1(V_1\cap V_2).$$ Now ${\cal
I}'={\cal I} -\{V_1,V_2\} + \{V_1\cap V_2\}$ is also $ST$-independent
and $\widetilde p_1({\cal I}')\geq \widetilde p_1({\cal I})$, but we
must have here equality by the optimality of $\cal I$, that is, ${\cal
I}'$ is also a $p_1$-minimizer, contradicting the minimality of $\vert
{\cal I}\vert $. $\bullet$\endproof

\begin{claim}\label{reduced.main} Let $\cal I$ be a strongly
$ST$-independent $p_1$-optimizer.  There exists a subset $X$ and a
subpartition ${\cal T}=\{T_1,\dots ,T_q\}$ of $T$ such that \ ${\cal
I}= \{V-s:  s\in X\} \cup \{X\cup T_i:  i=1,\dots ,q\}$, and hence \begin{equation}
\nu _1=\widetilde p_1({\cal I})= \widetilde m_S(X) + \widetilde
p_T({\cal T}) - \vert {\cal T}\vert \vert X\vert .
\label{(reduced.main)} \end{equation} \end{claim}

\proof{Proof.} Let $X:=\{s\in S:  V-s\in {\cal I}\}$ and let ${\cal
I}_1=\{V-s:  V-s\in {\cal I}\}$.  Let ${\cal I}_2:={\cal I} - {\cal
I}_1$ and let $V_1,\dots ,V_q$ denote the members of ${\cal I}_2$.
Furthermore, let $T_i:=T\cap V_i$ and $X_i = S\cap V_i$ \ ($i=1,\dots
,q)$.  By the strong $ST$-independence, the family ${\cal
T}=\{T_1,\dots ,T_q\}$ is a subpartition of $T$, and we also have
$X\subseteq X_i$ for each $i$.

Define $V_i':=T_i\cup X$ for $i=1,\dots ,q$ and let ${\cal
I}_2'=\{V_1',\dots ,V_q'\}$.  Then ${\cal I}'={\cal I}_1 \cup {\cal
I}_2'$ is also $ST$-independent.  Since $p_1(V_i')=p_1(V_i) +\vert
X_i-X\vert $ and $\cal I$ is a $p_1$-optimizer, we must have $X_i=X$
for each $i=1,\dots ,q$.  The formula in \eqref{(reduced.main)} follows
from $$ \nu _1= \widetilde p_1({\cal I}) = \widetilde p_1({\cal I}_1)
+ \widetilde p_1({\cal I}_2) =$$ $$ \sum [m_S(s):  V-s\in {\cal I}_1]
+ [\widetilde p_T({\cal T}) - \vert {\cal T}\vert \vert X\vert ] =
\widetilde m_S(X) + \widetilde p_T({\cal T}) - \vert {\cal T}\vert
\vert X\vert . \hbox{ $\bullet $ }\ $$\endproof

\begin{claim}\label{nugamma} $\nu _1=\gamma $. \end{claim}

\proof{Proof.} Since the family ${\cal L}=\{V-s:s\in S\}$ is $ST$-independent,
$\nu _1\geq \widetilde p_1({\cal L}) = \widetilde m_S(S) = \gamma $
from which $\nu _1\geq \gamma $. Let $\cal I$ be a strongly
$ST$-independent $p_1$-optimizer for which $\vert {\cal I}\vert $ is
minimum.  It follows from \eqref{(reduced.main)} in Claim
\ref{reduced.main} and from the hypothesis
\eqref{(basic.basefelt.1.free)} that $\nu _1\leq \gamma $ and hence
$\nu _1=\gamma $. $\bullet$\endproof \medskip

By Theorem \ref{Frank-Jordan}, there is a digraph $D=(V,A)$ on $V$
with $\nu _1=\gamma $ (possibly parallel) $ST$-arcs that covers $p_1$,
that is, $\varrho _D(V')\geq p_1(V')$ for every subset $V'\subseteq
V$.  Let $G=(S,T;E)$ denote the underlying bipartite graph of $D$.

\begin{claim}$d_G(s)=m_S(s)$ for every $s\in S$.\end{claim}

\proof{Proof.} Since $d_G(s)=\delta _D(s)= \varrho _D(V-s)\geq
p_1(V-s)=m_S(s)$ for every $s\in S$, we have $\gamma =\vert E\vert =
\sum [ d_G(s):s\in S] \geq \widetilde m_S(S) =\gamma $, from which
$d_G(s)=m_S(s)$ follows for every $s\in S$.  $\bullet$\endproof

\begin{claim}$\vert \Gamma _G(Y)\vert \geq p_T(Y)$ for every subset
$Y\subseteq T$.  \end{claim}

\proof{Proof.} Let $X:=\Gamma _G(Y)$ and $V':=X\cup Y$.  Then $0=\varrho
_D(V')\geq p_1(V')\geq p_0(V') = p_T(Y) - \vert X\vert = p_T(Y) -
\vert \Gamma _G(Y)\vert $, as required.  $\bullet$\endproof \medskip

Therefore the bipartite graph $G$ meets all the requirements of the
theorem apart possibly from simplicity.  By Claim \ref{simplify}, $G$
can be chosen to be simple.  $\bullet $ $\bullet $\endproof

\subsection{Degree-specification on \texorpdfstring{$S\cup T$}{S U T}} \label{ptst}

In the next problem we have degree-specification not only on $S$ but
on $T$ as well.  When the degree-specification was given only on $S$,
we have observed that it sufficed to concentrate on finding a
not-necessarily simple graph covering $p_T$ because such a graph could
easily be made simple.  Based on this, it is tempting to conjecture
that if there is a simple bipartite graph fitting a
degree-specification $m_V=(m_S,m_T)$ and there is a (not-necessarily
simple) one fitting $m_V$ and covering $p_T$, then there is a simple
bipartite graph fitting $m_V$ and covering $p_T$.  The following
example shows, however, that this statement fails to hold.

\begin{figure}[t]
\centering
\begin{subfigure}[b]{0.45\textwidth}
  \centering
  \includegraphics[width=0.65\linewidth]{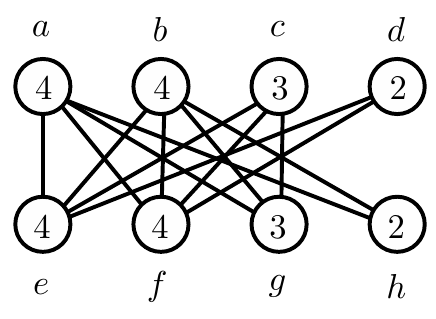}
  \caption{The unique simple graph fitting $m_V$}
  \label{fig:3.2}
\end{subfigure}%
\begin{subfigure}[b]{.45\textwidth}
  \centering
  \includegraphics[width=0.65\linewidth]{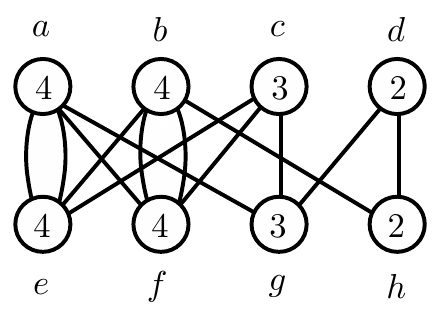}
  \caption{A non-simple graph covering $p_T$}
  \label{fig:3.2b}
\end{subfigure}
\caption{The existence of a simple graph fitting $m_V$ and a not-necessarily simple graph covering $p_T$ does not imply the existence of a graph satisfying these conditions simultaneously}
\label{fig:example}
\end{figure}

Let $S=\{e,f,g,h\}$ and let the $m_S$-values on $S$, respectively, be
$4,4,3,2.$ Let $T=\{a,b,c,d\}$ and let the $m_T$-values on $T$,
respectively, be $4,4,3,2.$ Let $p_T(t)=3$ for $t\in \{a,b,c\}$ and
let $p_T(d):=2$.  Let $p_T(\{c,d\})= 4$ and $p_T(\{y,z\})=1$ whenever
$\{y,z\}\not =\{c,d\}$, $\{y,z\}\subset T$.  Let
$p_T(\{a,c,d\})=p_T(\{b,c,d\}) = 3,$ and
$p_T(\{a,b,c\})=p_T(\{a,b,d\})=2$.  Finally, let $p_T(T)=4.$ Simple
case checking shows that function $p_T$ is $T$-intersecting
supermodular.  Here there is a unique simple bipartite graph $G$
fitting $m_V$, but $G$ does not cover $p_T$ since $\vert \Gamma
_G(\{c,d\})\vert = \vert \{e,f,g\}\vert =3\not \geq 4=p_T(\{c,d\})$.
On the other hand, the (non-simple) bipartite graph $G'=(S,T;E')$ with
$E'= \{ae,ae,af,ag, \ be, bf, bf, bh, \ ce, cf, cg, \ dg,dh\}$ fits
$m_V$ and covers $p_T$ (see Figure~\ref{fig:example}).

\begin{theorem}\label{ujterm.main.spec} Let $S$ and $T$ be disjoint sets and
let $m_V=(m_S,m_T)$ be a degree-specification for which $\widetilde
m_S(S)=\widetilde m_T(T)=\gamma $. Let $p_T$ be a positively
intersecting supermodular function on $T$ for which $p_T(\emptyset
)=0$ and $p_T(Y)\leq \vert S\vert $ for $Y\subseteq T$.  There is a
simple bigraph $G=(S,T;E)$ covering $p_T$ and fitting the
degree-specification $m_V$ if and only if 
\begin{eqnarray}
&\widetilde m_S(X) + \widetilde m_T(Y) -\vert X\vert \vert Y\vert + \widetilde p_T({\cal
T}) - \vert {\cal T}\vert \vert X\vert \leq \gamma & \nonumber\\[5pt]
& \hbox{ for every $X\subseteq S$, $Y\subseteq T$, and subpartition ${\cal T}$ of $T-Y$
} & \label{(ujterm.felt.1.spec)} 
\end{eqnarray} 
holds including the special case
${\cal T}=\emptyset $ {\em (when the condition is exactly
\eqref{(ujterm.GR)})}.  When $p_T$ is fully supermodular, it suffices
to require \eqref{(ujterm.felt.1.spec)} only for $\vert {\cal T}\vert
\leq 1$.  When $p_T$ is fully supermodular and monotone
non-decreasing, it suffices to require \eqref{(ujterm.felt.1.spec)}
only for ${\cal T}=\emptyset $ and for ${\cal T}=\{T-Y\}.$ \end{theorem}

\proof{Proof.} Necessity.  Suppose that there is a requested bigraph $G$.  Let
$X\subseteq S$, $Y\subseteq T$ be subsets and let ${\cal
T}=\{T_1,\dots ,T_q\}$ be a subpartition of $V-Y$.  On one hand, the
simplicity of $G$ implies that the number of edges with at least one
end-node in $X\cup Y$ is at least $\widetilde m_S(X) + \widetilde
m_T(Y) -\vert X\vert \vert Y\vert $. On the other hand, it was shown
already in \eqref{(pTdG)} that $d_G(T_i,S-X)\geq p_T(T_i) -\vert X\vert
$, implying that the number of edges between $S-X$ and $S-Y$ is at
least $\sum _i [ p_T(T_i) -\vert X\vert ]$.  Therefore the total
number $\gamma $ of edges of $G$ is at least $\widetilde m_S(X) +
\widetilde m_T(Y) -\vert X\vert \vert Y\vert + \widetilde p_T({\cal
T}) - \vert {\cal T}\vert \vert X\vert ,$ which is
\eqref{(ujterm.felt.1.spec)}.

Sufficiency.  Let $t$ be an element of $T$.  By applying
\eqref{(ujterm.felt.1.spec)} to $X=\emptyset , Y=T-t, q=1, T_1=\{t\}$,
we obtain that $\widetilde m_T(T-t) + p_T(t)\leq \gamma $, that is,
$p_T(t)\leq m_T(t)$.

Define a set-function $p_T\sp +$ on $T$ by revising $p_T$ so as to
lift its value on each singleton $\{t\}$ \ to $m_T(t)$ \ ($t\in T$).
As $p_T(t)\leq m_T(t)$ and $p_T$ is positively $T$-intersecting
supermodular, so is $p_T\sp +$.

Let $s$ be an element of $S$.  By applying \eqref{(ujterm.felt.1.spec)}
to $X=\{s\}$, $Y=T$, and $q=0$, we obtain that $m_S(s) + \widetilde
m_T(T) - \vert T\vert \leq \gamma $, that is, $m_S(s)\leq \vert T\vert
$, implying that \eqref{(basic.msb.free)} holds.

\begin{claim}Condition \eqref{(basic.basefelt.1.free)} holds for $p_T\sp +$
in place of $p_T$.\end{claim}

\proof{Proof.} Let $X\subseteq S$ and let ${\cal T}'=\{T_1,T_2,\dots
,T_{q'}\}$ be a subpartition of $T$.  Let $T_1,T_2,\dots ,T_q$ denote
those members of ${\cal T}'$ for which $p_T\sp +(T_i) = p_T(T_i)$ and
let ${\cal T}= \{T_1,T_2,\dots ,T_q\}$.  Then each of the remaining
members $T_j$ in ${\cal T}'$ is a singleton $\{z_j\}$ \ $(j=q+1,\cdots
,q')$ for which $p_T\sp +(T_j) = m_T(z_j)$.  By letting
$Y=\{z_{q+1},\dots ,z_{q'}\}$, we have $\vert Y\vert =q'-q$.  By
applying \eqref{(ujterm.felt.1.spec)} to this choice of $(X,Y,{\cal
T})$, we obtain that

$$\widetilde m_S(X) + \sum [p_T\sp +(T_i)-\vert X\vert :  i=1,\dots
,q'] =$$ $$\widetilde m_S(X) + \sum [p_T(T_i)-\vert X\vert :
i=1,\dots ,q] + \sum [m_T(z_j) - \vert X\vert :  j = q+1,\dots ,q']=$$
$$ \widetilde m_S(X) + \sum [p_T(T_i)-\vert X\vert :  i=1,\dots ,q] +
\widetilde m_T(Y) - \vert X\vert \vert Y\vert \leq \gamma ,$$

\noindent that is, condition \eqref{(basic.basefelt.1.free)} holds
indeed for $p_T\sp +$.

By applying Theorem \ref{basic.base.free} to $p_T\sp +$, we obtain
that there is a simple bipartite graph fitting the
degree-specification $m_S$ for which $\vert \Gamma _G(Y)\vert \geq
p_T\sp +(Y) \geq p_T(Y)$ for every subset $Y\subseteq T$.  In
particular, this implies for $Y=\{t\}$ that $d_G(t) = \vert \Gamma
_G(t)\vert \geq p_T\sp +(t) =m_T(t)$.  Therefore $\gamma = \sum
[d_G(t):t\in T] \geq \sum [m_T(t):t\in T] = \widetilde m_T(T) = \gamma
$ and hence we must have $d_G(t) = m_T(t)$ for every $t\in T$, making
the proof of the main part of the theorem complete.

\medskip Suppose now that $p_T$ is fully supermodular.  For specified
$X\subseteq S$ and $Y\subseteq T$, let $\cal T$ be a maximizer
subpartition in the left-hand side of the inequality in
\eqref{(ujterm.felt.1.spec)}.  We are done if $\vert {\cal T}\vert \leq
1$.  Suppose that ${\cal T}=\{T_1,\dots ,T_q\}$ for $q\geq 2$, and
consider the subpartition ${\cal T}'$ consisting of the single set
$T_0:=T_1\cup \cdots \cup T_q$.  The full supermodularity of $p_T$
implies $\widetilde p_T({\cal T})\leq \widetilde p_T({\cal T}')$.
Since $\cal T$ is a maximizer, we have $\widetilde p_T({\cal T})
-q\vert X\vert \geq \widetilde p_T({\cal T}') -\vert X\vert \geq
\widetilde p_T({\cal T}) -\vert X\vert \geq \widetilde p_T({\cal T})
-q\vert X\vert .$ Hence equality follows throughout, in particular,
$X=\emptyset $ and $\widetilde p_T({\cal T}) = \widetilde p_T({\cal
T}')$, showing that ${\cal T}'$ is also a maximizer.

Finally, investigate the case when $p_T$ is fully supermodular and
monotone non-decreasing.  If there are sets $X$, $Y$ and a
subpartition ${\cal T}$ of $T-Y$ violating \eqref{(ujterm.felt.1.spec)}
so that ${\cal T}=\{T_1\}$, then $X$, $Y$, and ${\cal T}'=\{T-Y\}$
also violates \eqref{(ujterm.felt.1.spec)} since $p_T(T-Y)\geq
p_T(T_1)$.  $\bullet $ $\bullet $ $\bullet $\endproof

\begin{corollary}\label{msmtfullp} Let $S,T, m_S,m_T, \gamma ,$ and $p_T$ be
the same as in Theorem \ref{ujterm.main.spec} and assume that $p_T$ is
non-decreasing and fully supermodular.  There is a simple bigraph
covering $p_T$ and fitting $(m_S,m_T)$ if and only if

\begin{equation}\widetilde m_S(X) + \widetilde m_T(Y) -\vert X\vert \vert Y\vert
\leq \gamma \ \hbox{ whenever $X\subseteq S$, $Y\subseteq T$}\
\label{(ujbasic.nemures1)} \end{equation} and

\begin{equation}\widetilde m_S(X) + \widetilde m_T(Y) -\vert X\vert \vert Y\vert +
p_T(T-Y) - \vert X\vert \leq \gamma \ \hbox{ whenever $X\subseteq S$,
$Y\subset T.$ }\ \label{(ujbasic.nemures2)} \end{equation} \end{corollary}

\proof{Proof.} Recall that the members of $\cal T$ in
\eqref{(ujterm.felt.1.spec)} are non-empty, in particular, if ${\cal
T}=\{T-Y\}$, then $Y\subset T$.  By the last part of Theorem
\ref{ujterm.main.spec}, the corollary follows.

\medskip

\begin{remark} In the example above, the subsets $X=\{e,f\}, Y=\{a,b\}$
and the subpartition ${\cal T}= \{ \{c,d\} \}$ consisting of a single
set (that is, $q=1$) do violate the necessary condition
\eqref{(ujterm.felt.1.spec)} since $ \widetilde m_S(X) + \widetilde
m_T(Y) -\vert X\vert \vert Y\vert + \sum _{i=1}\sp q [p_T(T_i)-\vert
X\vert ] = 8 + 8 - 4 + [4-2] =14\not \leq 13= 4+4+3+2= \gamma $. 
\end{remark}

\medskip

The essence of the next corollary of Theorem \ref{ujterm.main.spec} is
that it suffices to require \eqref{(ujterm.felt.1.spec)} only for
subsets $X\subseteq S$ with the $j$ largest $m_S$-values.  We leave
out the straightforward proof which consists of pointing out the
equivalence of \eqref{(ujterm.felt.cor)} and
\eqref{(ujterm.felt.1.spec)}.

\begin{corollary}\label{ujterm.mainfree} Let $S,T,p_T,$ and $m_V=(m_S,m_T)$
be the same as in Theorem \ref{ujterm.main.spec}.  There is a simple
bipartite graph $G=(S,T;E)$ covering $p_T$ and fitting $m_V$ if and
only if \begin{equation}\widetilde m_T(Y) + \sum _{i=1}\sp q p_T(T_i) \leq \sum
_{s\in S} \min \{m_S(s), \vert Y\vert +q\} \label{(ujterm.felt.cor)}
\end{equation} holds for every subset $Y\subseteq T$ and subpartition
$\{T_1,\dots ,T_q\}$ of $T-Y$ (including the special case when $q=0$
or $Y=\emptyset $).  \end{corollary}

\subsection{An {\bf NP}-complete extension}

One may be wondering if the synthesis problem solved in Theorem
\ref{ujterm.main.spec} could possibly be extended to the corresponding
subgraph problem.  That is, the problem is to characterize the
situation when the requested bigraph $G$ (covering $p_T$) is a
subgraph of an initial bipartite graph $G_0=(S,T;E_0)$.  However such
an extension is unlikely to exist since it includes {\bf NP}-complete
problems.

To see this, let $G_0=(S,T;E_0)$ be a bipartite graph in which $\vert
S\vert =\vert T\vert +1$.  Define $m_T$ to be identically 2 on $T$ and
$m_S$ to be identically 2 on $S$ apart from two specified nodes
$s_1,s_2\in S$ where $m_S(s_1) = m_S(s_2)=1$.  Define $p_T(Y)=\vert
Y\vert +1$ for each non-empty $Y\subseteq T$ and let $p_T(\emptyset
)=0$.  Clearly, $p_T$ is intersecting supermodular.

\begin{lemma}A subgraph $G=(S,T;E)$ of $G_0$ covers $p_T$ and fits
$m_V=(m_S,m_T)$ if and only if $G$ is a Hamiltonian path connecting
$s_1$ and $s_2$.\end{lemma}

\proof{Proof.} A Hamiltonian path $G$ contains a matching covering $T$ and
hence $\vert \Gamma _G(Y)\vert \geq \vert Y\vert $ for every
$Y\subseteq T$.  If indirectly $G$ does not cover $p_T$, then there is
a non-empty subset $Y$ of $T$ for which $\vert \Gamma _G(Y)\vert
=\vert Y\vert $. But then the subgraph of $G$ induced by $Y\cup \Gamma
_G(Y)$ has exactly $2\vert Y\vert = \vert Y\cup \Gamma _G(Y)\vert $
edges, contradicting the assumption that $G$ is a path.

Suppose now that $G$ covers $p_T$ and fits $m_V$.  Then $G$ has
$2\vert T\vert =\vert S\cup T\vert -1$ edges.  It cannot comprise a
circuit $C$ since then we would have $\vert \Gamma _G(Y)\vert = \vert
Y\vert $ for $Y=T\cap C$ contradicting the assumption that $G$ covers
$p_T$.  Therefore $G$ is a spanning tree, and since $G$ fits $m_V$, it
must be a Hamiltonian path connecting $s_1$ and $s_2$.  $\bullet$\endproof
\medskip

Since the Hamiltonian path problem is {\bf NP}-complete, so is the
equivalent problem of finding a subgraph of $G_0$ that covers $p_T$
and fits $m_V.$

Note that the same example shows that the synthesis problem solved in
Theorem \ref{basic.base.free} cannot be extended either to the
corresponding subgraph problem.

\medskip 

\section{The master base-polyhedron associated with realizable degree-specifications} \label{poli}

As before, $S$ and $T$ are two disjoint non-empty sets, $V:=S\cup T$,
and $m=(m_S,m_T)$ is a degree-specification for which $\widetilde
m_S(S) = \widetilde m_T(T)=\gamma $. Let $p_T$ be a positively
intersecting supermodular set-function on $T$ for which \begin{equation}\hbox{ $
p_T(Y)\leq \vert S\vert $ for every subset $Y\subseteq T$.}\
\label{(atmostS)} \end{equation} This implies that the complete bipartite graph
$(S,T;E\sp *)$ is a simple bigraph covering $p_T$.  Recall Theorem
\ref{ujterm.main.spec} which stated that there is a simple bigraph
covering $p_T$ and fitting $m$ if and only if \begin{eqnarray} & \widetilde
m_S(X) + \widetilde m_T(Y) - \vert X\vert \vert Y\vert + \widetilde
p_T({\cal T}) - \vert {\cal T}\vert \vert X\vert \leq \gamma & \nonumber \\[5pt]
& \hbox{whenever $X\subseteq S, \ Y\subseteq T$, and $\cal T$ a subpartition
of $T-Y$.} & \ \label{(mhezkell)} \end{eqnarray} We allow throughout the empty
subpartition with the convention $\widetilde p_T(\emptyset )=0$.  For
brevity we call such a degree-specification realizable (with respect
to $p_T$).  In this section, we investigate the problem when, rather
than an exact degree specification $m$, lower and upper bounds are
prescribed for the degrees of the requested simple bigraph covering
$p_T$.  Instead of attacking the problem directly, we exhibit first a
novel construction for a submodular function $b_0$ and show that there
is a simple one-to-one correspondence between the realizable
degree-specifications and the integral elements of the base-polyhedron
$B_0=B(b_0)$.  Because of its central role, we call $B_0$ the {\bf
master base-polyhedron} associated with $p_T$ and $S$.

Recall that for a submodular function $b$ with $b(V)$ finite, the
polyhedron $B(b):=\{x\in {\bf R}\sp V:  \widetilde x\leq b, \
\widetilde x(V)=b(V)\}$ is called a {\bf base-polyhedron}, and we
speak of a {\bf 0-base-polyhedron} if $b(V)=0$.  Given this
correspondence at hand, we can apply some known characterizations for
the non-emptiness of the intersection of a g-polymatroid with a box
and with a plank.  This approach enables us to treat situations when,
in addition to degree-constraints, upper and lower bounds for the
total number of edges can also be prescribed.

\subsection{A new submodular function} \label{base-polyhedron}

With each vector $m=(m_S,m_T)$, we associate the vector
$m'=(m_S,-m_T)$.  Note that the property $\widetilde m_S(S)=\widetilde
m_T(T)$ is equivalent to $\widetilde m'(V)=0$.  The condition
\eqref{(mhezkell)} for the realizability of $m$ is equivalent to the
following.  
\begin{eqnarray}
&\widetilde m'(X\cup Z) \leq \vert T-Z\vert
\vert X\vert -\widetilde p_T({\cal T}) + \vert {\cal T}\vert \vert
X\vert&\nonumber\\[5pt]
& \hbox{whenever $X\subseteq S, Z\subseteq T$, and $\cal T$ a
subpartition of $Z.$ }& \label{(m'kell)} 
\end{eqnarray}
Define a set-function
$b_0$ on $V$ as follows.  For $X\subseteq S$ and $Z\subseteq T$, let
\begin{equation}\hbox{ $ b_0(X\cup Z):  = \min \{ \vert T-Z\vert \vert X\vert
-\widetilde p_T({\cal T}) + \vert {\cal T}\vert \vert X\vert :  \ \cal
T$ a subpartition of $Z\}$.}\ \label{(bdef)} \end{equation} Clearly,
\eqref{(m'kell)} is equivalent to \begin{equation}\hbox{ $\widetilde m'(U)\leq
b_0(U)$ whenever $U\subseteq V$.}\ \label{(m'ekvi)} \end{equation}

\begin{claim}\label{Vnulla} $b_0(\emptyset )=0$ and $b_0(V)=0$.  \end{claim}

\proof{Proof.} When $Z=\emptyset $, a subpartition of $Z$ is also empty, and
hence $b_0(\emptyset )$ is indeed zero.

For $X=S$ and $Z=T$, we have $b_0(V)= \min \{ -\widetilde p_T({\cal
T}) + \vert {\cal T}\vert \vert S\vert :  {\cal T}$ a subpartition of
$T\}$.  By choosing $\cal T$ to be empty, we see that the minimum is
at most $0$.  On the other hand $-\widetilde p_T({\cal T}) + \vert
{\cal T}\vert \vert S\vert \geq 0$ holds for every subpartition $\cal
T$ of $T$ since \eqref{(atmostS)} implies that $\widetilde p_T({\cal
T})\leq \vert {\cal T}\vert \vert S\vert $. Therefore $b_0(V)=0$.
$\bullet$\endproof \medskip

\begin{theorem}$b_0$ is fully submodular.  \end{theorem}

\proof{Proof.} Let $V_1=X_1\cup Z_1$ and $V_2=X_2\cup Z_2$ be two subsets of
$V$ with $X_i\subseteq S$ and $Z_i\subseteq T$ \ $(i=1,2)$, and let
${\cal T}_i$ denote an optimizer subpartition of $Z_i$ in the
definition of $b_0(V_i)$.  That is, $$b_0(V_i) = \vert T-Z_i\vert
\vert X_i\vert -\widetilde p_T({\cal T}_i) + \vert {\cal T}_i\vert
\vert X_i\vert .$$

Let ${\cal F}_0$ denote the multi-union of ${\cal T}_1$ and ${\cal
T}_2$, that is, each member of ${\cal T}_1$ and ${\cal T}_2$ occurs in
${\cal F}_0$, and if $W$ is in both ${\cal T}_1$ and ${\cal T}_2$,
then two copies of $W$ occur in ${\cal F}_0$.  Hence $\vert {\cal
T}_1\vert + \vert {\cal T}_2\vert = \vert {\cal F}_0\vert $.

An {\bf uncrossing step} consists of taking two properly intersecting
members $A$ and $B$ of the current family for which $p_T(A)>0$ and
$p_T(B)>0$ and replacing them by $A\cup B$ and $A\cap B$.  (Note that
a set $A$ with $p(A)\leq 0$ never participates in an uncrossing step.)

The {\bf uncrossing procedure} starts with ${\cal F}_0$ and repeatedly
performs uncrossing steps.  It is known that the uncrossing procedure
is finite (as the number of sets does not change while the total sum
of the squares of cardinalities strictly increases).  Let ${\cal F}_0,
{\cal F}_1,{\cal F}_2,\dots , {\cal F}_q$ denote the subsequent
families, that is, ${\cal F}_{j+1}$ arises by applying the uncrossing
step to two members of ${\cal F}_{j}$ (which are properly intersecting
and have strictly positive $p$-values).

\begin{claim}\label{fedes-szam} Every family ${\cal F}_j$ \ $(j=0,\dots ,q)$
covers each element of $Z_1\cap Z_2$ at most twice, each element of
the symmetric difference $Z_1 \ominus Z_2$ at most once, and no
element outside $Z_1\cup Z_2$.  \end{claim}

\proof{Proof.} The property clearly holds for $j=0$ and it is maintained
throughout since an uncrossing step does not affect the number of sets
containing any given element of $T$.  $\bullet$\endproof

\begin{claim}\label{twocopies} If the family ${\cal F}_h$ for some
$h=0,\dots ,q$ contains two copies of a set $W$, then each family
${\cal F}_j$ \ $(j=0,\dots ,q)$ contains two copies of $W$.  In
particular, $W\in {\cal T}_1$ and $W\in {\cal T}_2$.  \end{claim}

\proof{Proof.} By induction, it suffices to show that both ${\cal F}_{h+1}$
and ${\cal F}_{h-1}$ contain two copies of $W$.

By Claim \ref{fedes-szam}, no member of ${\cal F}_h$ can intersect
properly $W$, and therefore both copies of $W$ belong to ${\cal
F}_{h+1}.$ Similarly, Claim \ref{fedes-szam} implies that both copies
of $W$ must be in ${\cal F}_{h-1}$ since if the second copy of $W$ in
${\cal F}_h$ arises as the intersection or the union of two properly
intersecting members $A$ and $B$ of ${\cal F}_{h-1}$, then the
elements of $A\cap B$ would belong to $A,B,$ and $W$.  $\bullet$\endproof

\begin{claim}\label{uncross.intersect} Let $W$ be a member of ${\cal
F}_{j+1}$ arising as the intersection of two properly intersecting
members $A$ and $B$ of ${\cal F}_j$, and let $Y$ be any member of
${\cal F}_{j+1}\cup \cdots \cup {\cal F}_q$ intersecting $W$.  Then
$W\subset Y$.\end{claim}

\proof{Proof.} We say that a pair of elements of $T$ is {\bf non-separated} by
a family of sets if no member of the family contains exactly one of
the two elements.  Clearly, if a pair is non-separated, then it
remains so after an uncrossing step.

By Claim \ref{twocopies}, $W$ does not occur in two copies and hence
$Y\not =W$.  By Claim \ref{fedes-szam}, any two elements of $A\cap B$
are non-separated by ${\cal F}_j$ and hence by each of ${\cal
F}_{j+1},\dots , {\cal F}_q$, as well.  Therefore, as $Y$ intersects
$W$, it must properly include $W$.  $\bullet$\endproof

\begin{claim}\label{uncross.union} Let $W$ be a member of ${\cal F}_{j+1}$
arising as the union of two properly intersecting members $A$ and $B$
of ${\cal F}_j$.  Then $W$ has a subset belonging to ${\cal T}_1$ and
$W$ has a subset belonging to ${\cal T}_2$.  \end{claim}

\proof{Proof.} Suppose the claim fails to hold and let $j$ be the smallest
index occurring in a counter-example.  If both $A$ and $B$ would
belong to ${\cal F}_0$, then one of them is in ${\cal T}_1$ while the
other one in ${\cal T}_2$, as these families are subpartitions.  But
in this case the pair $(W,j)$ would not be a counter-example.

Therefore at least one of $A$ and $B$, say $A$, is not in ${\cal
F}_0$.  By Claim \ref{uncross.intersect}, $A$ could not arise as an
intersection at an uncrossing step, that is, $A$ arose as the union of
two sets.  By the minimality of $j$, $A$ has a subset belonging to
${\cal T}_1$ and $A$ has a subset belonging to ${\cal T}_2$.  As $W$
is a superset of $A$, $W$ also has a subset belonging to ${\cal T}_1$
and a subset belonging to ${\cal T}_2$.  $\bullet$\endproof \medskip

Define $$ {\cal L}:= \{W\in {\cal F}_q:  \ p_T(W)>0\}.$$ Clearly,
$\cal L$ is laminar.  Let ${\cal P}_1$ consist of the minimal members
of $\cal L$ which are subsets of $Z_1\cap Z_2$, with the convention
that if two copies of a set $W\subseteq Z_1\cap Z_2$ belong to $\cal
L$, then one of them is placed in ${\cal P}_1$.  Let ${\cal P}_2$
consist of the members of $\cal L$ which are not in ${\cal P}_1$.

\begin{claim}\label{subpart} ${\cal P}_1$ is a subpartition of $Z_1\cap Z_2$
and ${\cal P}_2$ is a subpartition of $Z_1\cup Z_2$.  \end{claim}

\proof{Proof.} Since $\cal L$ is laminar, its minimal members are disjoint and
hence ${\cal P}_1$ is indeed a subpartition.

To see that ${\cal P}_2$ is also a subpartition, assume indirectly
that two members $A$ and $B$ of ${\cal P}_2$ are not disjoint.  Then
the laminarity of $\cal L$ implies that one of $A$ and $B$ includes
the other, say, $A\subseteq B$.  We must have $A\subset B$ for if we
had $A=B$, then $A\subseteq Z_1\cap Z_2$ by Claim \ref{twocopies} and
one of $A$ and $B$ would belong to ${\cal P}_1$ by the definition of
${\cal P}_1$.  Because each element of $Z_1 \ominus Z_2$ belongs to at
most one member of $\cal L$, we have $A\subseteq Z_1\cap Z_2$.  But
$A$ is not in ${\cal P}_1$, that is, $A$ is not a minimal member of
$\cal L$, contradicting the property that each element of $T$ belongs
to at most two members of $\cal L$.  $\bullet$\endproof

\begin{claim}\label{T1T2} Let $W$ be a member of ${\cal P}_2$.  If
$W\subseteq Z_i$ \ $(i=1,2)$, then $W$ has a subset belonging to
${\cal T}_i$.  \end{claim}

\proof{Proof.} Since the indices 1 and 2 play a symmetric role, we prove the
claim only for $i=1$.  That is, we assume that $W\subseteq Z_1$ and
will show that there is a subset of $W$ belonging to ${\cal T}_1$.  If
$W$ is in ${\cal P}_1$, as well, that is, if two copies of $W$ occur
in $\cal L$, then we are done by Claim \ref{twocopies}.  Therefore, we
can assume that $W\not \in {\cal P}_1$.

By Claim \ref{uncross.union}, we are done if $W$ has arisen as a union
during the uncrossing procedure.  Suppose now that $W$ arises as an
intersection of $A$ and $B$ during the uncrossing procedure.  Then
Claim \ref{fedes-szam} implies that $W=A\cap B\subseteq Z_1\cap Z_2$.
Since $W$ is not in ${\cal P}_1$, there must be a set $Y\in {\cal L}$
for which $Y\subset W$, contradicting Claim \ref{uncross.intersect}.

In the remaining case, $W$ belongs each of the families ${\cal
F}_0,{\cal F}_1,\dots ,{\cal F}_q$.  In particular, $W$ is in ${\cal
F}_0$.  Since we are done if $W\in {\cal T}_1$, we can assume that
$W\in {\cal T}_2$.  In this case, $W-Z_2=\emptyset $, that is,
$W\subseteq Z_1\cap Z_2$.  Since $W$ is not in ${\cal P}_1$, there
must be a set $Y\in {\cal L}$ for which $Y\subset W$.  Since $W$
belongs to each ${\cal F}_j$, $Y$ could not arise as an intersection
or a union during the uncrossing procedure, and therefore $Y$ is also
a member of ${\cal F}_0$.  Since ${\cal T}_2$ is a subpartition, $Y$
cannot be in ${\cal T}_2$, that is, $Y\in {\cal T}_1$.  $\bullet$\endproof
\medskip

For simplifying calculations, we introduce the following four
parameters.

$$ \hbox{ $\tau _1:= \vert T-Z_1\vert + \vert {\cal T}_1\vert $ \ and
\ $\tau _2:= \vert T-Z_2\vert + \vert {\cal T}_2\vert ,$}\ $$ $$
\hbox{ $\pi _1:= \vert T-(Z_1\cap Z_2)\vert + \vert {\cal P}_1\vert $
\ and \ $\pi _2:= \vert T-(Z_1\cup Z_2)\vert +\vert {\cal P}_2\vert
.$}\ $$

\begin{claim}\label{pi2korlat} $\pi _2\leq \tau _1$ \ and \ $\pi _2\leq \tau
_2$.  \end{claim}

\proof{Proof.} Since the role of $\tau _1$ and $\tau _2$ is symmetric, we
prove only the first inequality.  Since ${\cal P}_2$ is a
subpartition, ${\cal P}_2$ has at most $\vert Z_2-Z_1\vert $ members
intersecting $Z_2-Z_1$, and, by Claim \ref{T1T2}, ${\cal P}_2$ has at
most ${\cal T}_1$ members not intersecting $Z_2-Z_1$.  Therefore
$\vert {\cal P}_2\vert \leq \vert Z_2-Z_1\vert + \vert {\cal T}_1\vert
$. By adding this to the identity $\vert T-(Z_1\cup Z_2)\vert = \vert
T-Z_1\vert - \vert Z_2-Z_1\vert $, we obtain the required $\pi _2\leq
\tau _1.$ $\bullet$\endproof

\begin{claim}\label{TiPi} $$\tau _1+\tau _2 \geq \pi _1+\pi _2$$ and
$$\widetilde p_T({\cal T}_1)+ \widetilde p_T({\cal T}_2) \leq
\widetilde p_T({\cal P}_1)+ \widetilde p_T({\cal P}_2).$$ \end{claim}

\proof{Proof.} Clearly, $\vert {\cal T}_1\vert + \vert {\cal T}_2\vert = \vert
{\cal F}_0\vert = \vert {\cal F}_q\vert \geq \vert {\cal L}\vert =
\vert {\cal P}_1\vert + \vert {\cal P}_2\vert .$ By adding this to
$\vert T-Z_1\vert + \vert T-Z_2\vert = \vert T - (Z_1\cap Z_2)\vert +
\vert T - (Z_1\cup Z_2)\vert $, the first inequality follows.

Since $p_T$ is positively intersecting supermodular, an uncrossing
step cannot decrease the $p_T$-sum of the current family.  Hence
$\widetilde p_T({\cal T}_1)+ \widetilde p_T({\cal T}_2) = \widetilde
p_T({\cal F}_0) \leq \widetilde p_T({\cal F}_q) \leq \widetilde
p_T({\cal L})= \widetilde p_T({\cal P}_1)+ \widetilde p_T({\cal
P}_2).$ $\bullet$\endproof \medskip

For $i=1,2$, we have:  \begin{equation}b_0(V_i) = \vert T-Z_i\vert \vert X_i\vert
-\widetilde p_T({\cal T}_i) + \vert {\cal T}_i\vert \vert X_i\vert \ =
\ \tau _i\vert X_i\vert -\widetilde p_T({\cal T}_i).  \label{(bVi)}
\end{equation} Since ${\cal P}_1$ is a subpartition of $Z_1\cap Z_2$, we have
\begin{equation}b_0(V_1\cap V_2) \leq \vert T-(Z_1\cap Z_2)\vert \vert X_1\cap
X_2\vert -\widetilde p_T({\cal P}_1) + \vert {\cal P}_1\vert \vert
X_1\cap X_2\vert \ = \ \pi _1\vert X_1\cap X_2\vert -\widetilde
p_T({\cal P}_1).  \label{(bmetszet)} \end{equation} Since ${\cal P}_2$ is a
subpartition of $Z_1\cup Z_2$, we have \begin{equation}b_0(V_1\cup V_2) \leq \vert
T-(Z_1\cup Z_2)\vert \vert X_1\cup X_2\vert -\widetilde p_T({\cal
P}_2) + \vert {\cal P}_2\vert \vert X_1\cup X_2\vert \ = \ \pi _2\vert
X_1\cup X_2\vert -\widetilde p_T({\cal P}_2).  \label{(bunio)} \end{equation} By
combining these inequalities, we obtain:  
\begin{eqnarray*}
&b_0(V_1) + b_0(V_2) =
[\tau _1\vert X_1\vert -\widetilde p_T({\cal T}_1)] + [\tau _2\vert
X_2\vert -\widetilde p_T({\cal T}_2)] =&\\[8pt] &\tau _1\vert
X_1-X_2\vert + \tau _2\vert X_2-X_1\vert + (\tau _1+\tau _2)\vert
X_1\cap X_2\vert -\widetilde p_T({\cal T}_1) -\widetilde p_T({\cal
T}_2) \geq &\\[8pt] &\pi _2\vert X_1-X_2\vert + \pi _2\vert X_2-X_1\vert +
(\pi _1+\pi _2)\vert X_1\cap X_2\vert -\widetilde p_T({\cal P}_1)
-\widetilde p_T({\cal P}_2) =& \\[8pt] & [\pi _1\vert X_1\cap X_2\vert
-\widetilde p_T({\cal P}_1)] + [\pi _2\vert X_1\cup X_2\vert
-\widetilde p_T({\cal P}_2)] \geq & \\[8pt] & b_0(V_1\cap V_2) + b_0(V_1\cup
V_2),
\end{eqnarray*}
that is, the function $b_0$ is indeed fully submodular.
$\bullet $ $\bullet $ \endproof\medskip

\begin{corollary}\label{mB-ben} An integral vector $m=(m_S,m_T)$ is the
degree-vector of a simple bigraph covering $p_T$ if and only if the
associated vector $m'=(m_S,-m_T)$ belongs to the 0-base-polyhedron
$B(b_0):=\{x\in {\bf R}\sp V:  \ \widetilde x\leq b_0, \ \widetilde
x(V)=0\}$. \end{corollary}

The following claim appeared in \cite{FrankP6} (see also Theorem
14.2.2 in book \cite{Frank-book}).

\begin{claim}\label{project} Given a non-empty subset $S\subset V,$ the
projection $Q'$ of a g-polymatroid $Q=Q(p,b)$ to ${\bf R}\sp S$ (or,
for short, to $S$) is the g-polymatroid $Q(p\vert _S,b\vert _S)$ where
$p\vert _S$ and $b\vert _S$ are the restriction of $p$ and $b$,
respectively, on $S$.  Each integral element of $Q'$ is the projection
of an integral element of $Q$.  \end{claim}

\begin{corollary}\label{csak.S-en} There is an integral g-polymatroid $Q_S$
in ${\bf R}\sp S$ so that a vector $m_S:S\rightarrow {\bf Z}_+$
belongs to $Q_S$ if and only if there is a simple bigraph covering
$p_T$ for which $d_G(s)=m_S(s)$ for every $s\in S$.  \end{corollary}

\proof{Proof.} Take $Q_S$ to be the projection of $B(b_0)$ to $S$ and apply
Claim \ref{project}.  $\bullet$\endproof

\section{Degree and edge-number constraints}

\subsection{Basic properties of generalized polymatroids}

In what follows, we make use of some basic notions and theorems of the
theory of generalized polymatroids.  (For a background, see for
example \cite{FrankJ17} or Chapter 14 in book \cite{Frank-book}.)  Let
$(p,b)$ be a fully paramodular (or, for short, paramodular) pair of
set-functions $p$ and $b$ defined on a ground-set $V$.  By definition,
this means that $b$ is submodular, $p$ is supermodular, and $$b(X)-p(Y) \geq b(X-Y) - p(Y-X)$$ holds for every
pair of subsets $X,\ Y$ of $V$.  The polyhedron $Q(p,b):=\{ x\in {\bf
R}\sp V:  p\leq \widetilde x\leq b \}$ is called a g-polymatroid and
$(p,b)$ is its {\bf border pair}.  Here we consider only
integer-valued functions $p$ and $b$.  The empty set is also
considered as a g-polymatroid, though it cannot be defined with the
help of a paramodular pair.  A special g-polymatroid is a box $T(f,g)
=\{x\in {\bf R}\sp V:  f\leq x\leq g\}$ where $f:V\rightarrow {\bf
Z}\cup \{-\infty \}$, $g:V\rightarrow {\bf Z}\cup \{\infty \}$ with
$f\leq g$.  Another special g-polymatroid is a plank $K(\alpha ,\beta
) =\{x\in {\bf R}\sp V:  \alpha \leq \widetilde x(V)\leq \beta \}$
where $\alpha \in {\bf Z}\cup \{-\infty \}$, $\beta \in {\bf Z}\cup
\{+\infty \}$ with $\alpha \leq \beta $.

With a submodular function $b$ with finite $b(V)$, we can associate
its {\bf complementary set-function} $p$ defined for $U\subseteq V$ by
$p(U):=b(V)-b(V-U).$ We list some basic properties.

\begin{claim}If $p$ is the complementary function of a submodular function
$b$, then $(p,b)$ is paramodular and $B(b)=Q(p,b)$.\end{claim}

\begin{claim}\label{egeszgpol} A g-polymatroid defined by an integral
paramodular pair is a non-empty integral polyhedron.  \end{claim}

\begin{claim}\label{pubu} A non-empty $g$-polymatroid $Q$ uniquely
determines its defining paramodular pair $(p,b)$, namely, $$ \hbox{
$p(U) = \min \{\widetilde x(U):  x\in Q\}$ and $b(U) = \max
\{\widetilde x(U):  x\in Q\}.$ }\ $$ \end{claim} \begin{claim}\label{metszet} The
intersection of two integral g-polymatroids is an integral polyhedron.
$Q(p_1,b_1)\cap Q(p_2,b_2)$ is non-empty if and only if $p_1\leq b_2$
and $p_2\leq b_1$.  \end{claim}

\begin{claim}The intersection of a g-polymatroid, a box, and a plank is a
g-polymatroid.\end{claim}

\begin{claim}\label{boxmetszet} The intersection $Q'$ of a g-polymatroid
$Q=Q(p,b)$ and a box $T=T(f,g)$ is non-empty if and only if
$\widetilde f\leq b$ and $p\leq \widetilde g$.  When $Q'$ is
non-empty, its unique border pair $(p',b')$ is given by \begin{equation}p'(U) =
\max\{ p(U') - \widetilde {g}(U'-U)+ \widetilde {f}(U-U') :
U'\subseteq V\}, \label{(maxp)} \end{equation} \begin{equation}b'(U) = \min\{ b(U') -
\widetilde {f}(U'-U)+ \widetilde {g}(U-U') :  U'\subseteq V\}.
\label{(minb)} \end{equation} \end{claim}

\begin{claim}[Linking property of g-polymatroids] If a g-polymatroid
$Q=Q(p,b)$ has an element $x'$ with $x'\geq f$, and $Q$ has an element
$x''$ with $x''\leq g$, then $Q$ has an element $x$ with $f\leq x\leq
g$.  In addition, $x$ can be chosen to be integral if $p,b,f,g$ are
integral.\end{claim}

\begin{claim}\label{savmetszet} The intersection $Q'$ of g-polymatroid
$Q=Q(p,b)$ and a plank $K(\alpha ,\beta )$ is non-empty if and only if
$\alpha \leq b(S)$ and $p(S)\leq \beta $. In particular, if $Q$ has an
element $x'$ with $\widetilde x'(V)\geq \alpha $ and $Q$ has an
element $x''$ with $\widetilde x''(V)\leq \beta $, then $Q$ has an
element $x$ with $\alpha \leq \widetilde {x}(V)\leq \beta $. Moreover,
if $p,b,\alpha ,\beta $ are integral, then $Q'$ is an integral
polyhedron.  \end{claim}

\subsection{Degree constraints}

We are given a lower bound function $f_V=(f_S,f_T)$ and an upper bound
function $g_V=(g_S,g_T)$ on $V=S\cup T$ for which $-\infty \leq f_V
\leq g_V \leq +\infty $.

\begin{theorem}\label{fg.main} Let $p_T$ be a positively intersecting
supermodular function on $T$ for which $p_T(Y)\leq \vert S\vert $ for
every $Y\subseteq T$.  There is a simple bigraph $G=(S,T;E)$ covering
$p_T$ and degree-constrained by $(f,g)$ if and only if
\begin{eqnarray}
&\widetilde f_T(Y) - \vert X\vert \vert Y\vert + \widetilde p_T({\cal
T}) - \vert {\cal T}\vert \vert X\vert \leq \widetilde g_S(S-X)& \nonumber \\[5pt]
&\text{for
every $Y\subseteq T$, $X\subseteq S$, \ $\cal T$ a subpartition of
$T-Y$ }& \label{(ftgs)} \end{eqnarray} and \begin{eqnarray}&\widetilde f_S(X) -
\vert X\vert \vert Y\vert + \widetilde p_T({\cal T}) - \vert {\cal
T}\vert \vert X\vert \leq \widetilde g_T(T-Y) & \nonumber \\[5pt]
&\text{for every $Y\subseteq
T$, $X\subseteq S$, \ $\cal T$ a subpartition of $T-Y$.  }&
\label{(fsgt)} \end{eqnarray} If $p_T$ is fully supermodular, then it suffices
to require the two conditions only for subpartitions $\cal T$ having
at most one member.  If $p_T$ is fully supermodular and monotone
non-decreasing, then it suffices to require the two conditions only
for ${\cal T}= \{\emptyset \}$ and ${\cal T}=\{T-Y\}.$ \end{theorem}

\proof{Proof.} Let \begin{equation}\hbox{ $f':=(f_S,-g_T)$ and $g':=(g_S,-f_T)$.}\
\label{(f'g')} \end{equation} Recall the submodular function $b_0$ and let $p_0$
denote its complementary function (that is, $p_0(U)=-b_0(V-U))$.  Then
$B(b_0)=Q(p_0,b_0)$ and, by Corollary \ref{mB-ben}, the requested
bigraph exists if and only if the intersection $Q'=Q(p_0,b_0)\cap
T(f',g')$ is non-empty.  By Claim \ref{boxmetszet}, $Q'$ is non-empty
precisely if $\widetilde f'\leq b_0$ and $p_0\leq \widetilde {g'}$.
We are going to show that $\widetilde f'\leq b_0$ is equivalent to
\eqref{(fsgt)} and that $p_0\leq \widetilde {g'}$ is equivalent to
\eqref{(ftgs)}.

By \eqref{(bdef)}, $\widetilde f'\leq b_0$ is equivalent to requiring
the following inequality for every pair of subsets $X\subseteq S, \
Z\subseteq T$:

$$ \hbox{ $\widetilde f'(X\cup Z) \leq \vert T-Z\vert \vert X\vert -
\widetilde p_T({\cal T}) + \vert {\cal T}\vert \vert X\vert $ whenever
$\cal T$ is a subpartition of $Z$.  }\ $$

By taking $Y:=T-Z$ and observing that $\widetilde f'(X\cup Z)=
\widetilde f_S(X) - \widetilde g_T(Z) = \widetilde f_S(X) - \widetilde
g_T(T-Y)$, we conclude that $\widetilde f'\leq b_0$ is equivalent to
\begin{eqnarray*}
&\widetilde f_S(X) - \widetilde g_T(T-Y) \leq \vert Y\vert
\vert X\vert - \widetilde p_T({\cal T}) + \vert {\cal T}\vert \vert
X\vert&\nonumber\\[5pt]
& \hbox{whenever $X\subseteq S, Y\subseteq T$, and $\cal T$ a subpartition
of $T-Y$,} & 
\end{eqnarray*}
which is the same as \eqref{(fsgt)}.

Let us prove now the equivalence of $p_0\leq \widetilde {g'}$ and
\eqref{(ftgs)}.  By taking $Y:=T-Z$ and $X':=S-X$, we have $\widetilde
{g'}(X'\cup Y)= \widetilde g_S(X')-\widetilde f_T(Y)= \widetilde
g_S(S-X) - \widetilde f_T(Y)$ and $$ \hbox{ $p_0(X'\cup Y) =
-b_0(X\cup Z) = - \min \{ \vert T-Z\vert \vert X\vert -\widetilde
p_T({\cal T}) + \vert {\cal T}\vert \vert X\vert :  \cal T$ a
subpartition of $Z\}.$ }\ $$

Condition $\widetilde g'\geq p_0$ means that $\widetilde g'(X'\cup Y)
\geq p_0(X'\cup Y)$ for every pair of sets $X'\subseteq S, Y\subseteq
T$, and this is equivalent to requiring $$ \widetilde g_S(S-X) -
\widetilde f_T(Y) \geq - [ \vert Y\vert \vert X\vert -\widetilde
p_T({\cal T}) + \vert {\cal T}\vert \vert X\vert ]$$ for every
subpartition $\cal T$ of $T-Y$, and this inequality is the same as the
one in \eqref{(ftgs)}.

The last part of the theorem concerning fully supermodular $p_T$
follows exactly the same way how the analogous statement was derived
in the proof of Theorem \ref{ujterm.main.spec}.  $\bullet$\endproof

\begin{corollary}\label{bigraph-linking} Let $p_T$ be a positively
intersecting supermodular function on $T$ for which $p_T(Y)\leq \vert
S\vert $ whenever $Y\subseteq T$.

\medskip \noindent {\bf (A)} \ There is a simple bigraph $G'$ covering
$p_T$ and degree-constrained by $(f_T,g_S)$ if and only if
\eqref{(ftgs)} holds.  \medskip

\noindent {\bf (B)} \ There is a simple bigraph $G''$ covering $p_T$
and degree-constrained by $(f_S,g_T)$ if and only if \eqref{(fsgt)}
holds.  \medskip

\noindent {\bf (AB)} There is a simple bigraph $G$ covering $p_T$ and
degree-constrained by $(f_V,g_V)$ if and only if both $G'$ and $G''$
exist (that is, both \eqref{(ftgs)} and \eqref{(fsgt)} hold).  \medskip

When $p_T$ is fully supermodular, it suffices to require the two
conditions only for subpartitions $\cal T$ having at most one member.
If $p_T$ is fully supermodular and monotone non-decreasing, then it
suffices to require the two conditions only for ${\cal T}= \{\emptyset
\}$ and ${\cal T}=\{T-Y\}.$ \end{corollary}

\proof{Proof.} (A) Define $f_S:\equiv -\infty $ and $g_T:\equiv +\infty $, and
observe that \eqref{(fsgt)} automatically holds when $X\not =\emptyset
$ or $Y\subset T$.  If $X=\emptyset $ and $Y=T$, then $\cal T$ is
empty and the requirement in \eqref{(fsgt)} becomes void.  Hence
Theorem \ref{fg.main} implies Part (A).

(B) Define $f_T:\equiv -\infty $ and $g_S:\equiv +\infty $, and
observe that \eqref{(ftgs)} automatically holds when $X\subset S$ or
$Y\not =\emptyset $. If $X=S$ and $Y=\emptyset $, then \eqref{(fsgt)}
reduces to $\widetilde p_T{(\cal T)} \leq \vert {\cal T}\vert \vert
S\vert $ for every subpartition $\cal T$ of $T$, but this follows from
the hypothesis that $p_T(Y)\leq \vert S\vert $ for every $Y\subseteq
T$.  Hence Theorem \ref{fg.main} implies Part (B).

(AB) \ Theorem \ref{fg.main} implies immediately Part (AB).  $\bullet
$

\begin{corollary}\label{mandg} Let $S$ and $T$ be disjoint sets and let
$m_S$ be a degree-specification on $S$ for which $\widetilde m_S(S)=
\gamma $. Let $g_T:T\rightarrow {\bf Z}_+$ be an upper bound function
for which $g_T(t)\leq \vert S\vert $ for every $t\in T$.  Let $p_T$ be
a positively intersecting supermodular function on $T$ with
$p_T(\emptyset )=0$.  There is a simple bigraph covering $p_T$ and
fitting $m_S$ for which \begin{equation}d_G(t)\leq g_T(t) \ \hbox{whenever}\ t\in
T \label{(felso.gt)} \end{equation} if and only if \begin{equation}\hbox{ $\widetilde m_S(X)
+ \widetilde p_T({\cal T}) - \vert {\cal T}\vert \vert X\vert \leq
\gamma $ \ whenever $X\subseteq S$ and ${\cal T}$ a subpartition of
$T$ }\ \label{(felso.felt.1)} \end{equation} and \begin{eqnarray}&\widetilde m_S(X)
-\vert X\vert \vert Y\vert + \widetilde p_T({\cal T}) - \vert {\cal
T}\vert \vert X\vert \leq \widetilde g_T(T-Y)& \nonumber\\[5pt] &\text{whenever $X\subseteq
S, \ Y\subseteq T$, and $\cal T$ a subpartition of $T-Y$, }&
\label{(felso.felt.2)} \end{eqnarray} where each of $X$, $Y$, and $\cal T$ may
be empty.  \end{corollary}

\proof{Proof.} (outline) Define $f_S:=m_S$ $g_S:=m_S$, $f_T:\equiv -\infty $,
and apply Theorem \ref{fg.main}.  $\bullet$\endproof

\medskip Note that in Corollary \ref{mandg} there is no need to
require explicitly the necessary condition given in \eqref{(atmostS)}
since this is implied by applying \eqref{(felso.felt.1)} in the special
case $X:=S$ and ${\cal T}:= \{Y\}$.

\begin{corollary}\label{csak.fsgs} Let $p_T$ be a positively intersecting
supermodular function on $T$ for which $p_T(Y)\leq \vert S\vert $ whenever
$Y\subseteq T$.  There is a simple bigraph $G=(S,T;E)$ covering $p_T$
and degree-constrained by $(f_S,g_S)$ if and only if \begin{equation}\hbox{
$f_S(s)\leq \vert T\vert $ whenever $s\in S$ }\ \label{(fsT)} \end{equation} and \begin{equation}
\hbox{ $ \widetilde p_T({\cal T}) - \vert {\cal T}\vert \vert X\vert
\leq \widetilde g_S(S-X)$ \ whenever $X\subseteq S$ and $\cal T$ a
subpartition of $T$. }\ \label{(csakgs)} \end{equation} \end{corollary}

\proof{Proof.} Define $f_T(t):\equiv -\infty $ and $g_T(t):\equiv +\infty $
and apply Theorem \ref{fg.main}.  Observe that condition \eqref{(ftgs)}
automatically holds when $Y\not =\emptyset $. When $Y=\emptyset $,
\eqref{(ftgs)} is just \eqref{(csakgs)}.  Similarly, condition
\eqref{(fsgt)} automatically holds when $Y\not =T$.  When $Y=T$, then
${\cal T}=\emptyset $ and \eqref{(fsgt)} requires $f_S(X) \leq \vert
X\vert \vert T\vert $ for every $X\subseteq S$ but this is equivalent
to \eqref{(fsT)}.  $\bullet$\endproof

\begin{corollary}\label{csak.ftgt} Let $p_T$ be a positively intersecting
supermodular function on $T$ for which $p_T(Y)\leq \vert S\vert $ for
every $Y\subseteq T$.  Let $g_T:T\rightarrow {\bf Z}_+$ be a function
for which $g_T(t)\leq \vert S\vert $ for every $t\in T$.  There is a
simple bigraph $G=(S,T;E)$ covering $p_T$ and degree-constrained by
$g_T$ if and only if \begin{equation}\hbox{ $ p_T(Y) \leq \widetilde g_T(Y)$ for
every $Y\subseteq T$.  }\ \label{(csak.ftgt)} \end{equation} \end{corollary}

\proof{Proof.} Define $f_V:\equiv -\infty $ and $g_S:\equiv +\infty $. Then
\eqref{(ftgs)} holds automatically (as we showed this in the proof of
Part (B) of Corollary \ref{bigraph-linking}).  Condition \eqref{(fsgt)}
holds automatically when $X\not =\emptyset $. If $X= \emptyset $, then
\eqref{(fsgt)} transforms to $$ \hbox{ $\widetilde p_T({\cal T)}\leq
\widetilde g_T(T-Y)$ whenever $Y\subset T$ and ${\cal T}=\{V_1,\dots
,V_q\}$ a subpartition of $T-Y$.  }\ $$ By Condition
\eqref{(csak.ftgt)}, $p_T(V_i) \leq \widetilde g(V_i)$ from which
$\widetilde p_T({\cal T})\leq \sum _{i=1}\sp q\widetilde g_T(V_i) \leq
\widetilde g_T(T-Y)$.  Therefore the conditions of Theorem
\ref{fg.main} hold and hence the required degree-constrained bigraph
exists.  $\bullet$\endproof

\medskip

\begin{remark} Corollary \ref{csak.ftgt} is not particularly exciting
since it can actually be formulated in a more general form when $G$ is
a subgraph of an initial bipartite graph $G_0$.  That was the content
of Theorem \ref{Frank-Tardos}.  To derive Corollary \ref{csak.ftgt},
choose $G_0$ to be the complete bipartite graph $G\sp *=(S,T,E\sp *)$
and observe that \eqref{(felso.Fra-Tar.felt1)} holds automatically when
$Z\not =\emptyset $. For $Z=\emptyset $, \eqref{(felso.Fra-Tar.felt1)}
is just \eqref{(csak.ftgt)}. \end{remark}

\medskip

\subsection{Edge-number constraints }

Suppose now that there exists a simple bigraph covering $p_T$ and
constrained by $(f,g)$, that is, conditions \eqref{(ftgs)} and
\eqref{(fsgt)} hold.  Our next goal is to characterize the situation
when, in addition to the degree constraints $(f,g)$, there are lower
and upper bounds $\alpha \leq \beta $ for the number of edges, as
well, where $\alpha $ and $\beta $ are non-negative integers.

\begin{theorem}\label{fg.ab} Suppose that conditions \eqref{(ftgs)} and
\eqref{(fsgt)} hold.  There is simple bigraph $G=(S,T;E)$ covering
$p_T$ and degree-constrained by $(f,g)$ for which \medskip

\noindent {\bf (A)} \ $\alpha \leq \vert E\vert $ if and only if \begin{equation}
\left\{ \begin{array}{ll} \widetilde g_S(S-X) + \widetilde g_T(T-Y) +
\vert X\vert \vert Y\vert - [\widetilde p_T({\cal T}) - \vert X\vert
\vert {\cal T}\vert ] \geq \alpha \\ \hbox{whenever $X\subseteq S, \
Y\subseteq T,$ and \ ${\cal T}$ a subpartition of $T-Y$,}\ &
\end{array} \right.  \label{(galfa)} \end{equation}

\noindent {\bf (B)} \ $\vert E\vert \leq \beta $ if and only if \begin{equation}
\left\{ \begin{array}{ll} \widetilde f_S(X) + \widetilde f_T(Y) -
\vert X\vert \vert Y\vert + \widetilde p_T({\cal T}) - \vert X\vert
\vert {\cal T}\vert \leq \beta \\ \hbox{whenever $X\subseteq S, \
Y\subseteq T$, and ${\cal T}$ a subpartition of $T-Y$,}\ & \end{array}
\right.  \label{(fbeta)} \end{equation}

\noindent {\bf (AB)} \ $\alpha \leq \vert E\vert \leq \beta $ if and
only if both \eqref{(galfa)} and \eqref{(fbeta)} hold.  \medskip

When $p_T$ is fully supermodular, it suffices to require the two
conditions only for subpartitions $\cal T$ having at most one member.
If $p_T$ is fully supermodular and monotone non-decreasing, then it
suffices to require the two conditions only for ${\cal T}= \{\emptyset
\}$ and ${\cal T}=\{T-Y\}.$ \end{theorem}

\proof{Proof.} Consider the functions $f'$ and $g'$ defined in \eqref{(f'g')}.
As we proved above, there is a simple bigraph covering $p_T$ and
constrained by $(f,g)$ if and only if the g-polymatroid
$Q'=Q(p_0,b_0)\cap T(f',g')$ is non-empty.  By our hypothesis $Q'$ is
non-empty.  Let $(p',b')$ denote the unique border pair of $Q'$ which
can be obtained by applying Claim \ref{boxmetszet} to $p_0,b_0,
f',g'$.

Let $Q'_S$ denote the projection of $Q'$ to $S$.  By Claim
\ref{project} the unique border pair of $Q'_S$ is $(p'\vert _S,b'\vert
_S)$, and any integral element of $Q'_S$ is the projection of an
integral element of $Q'$.  Therefore the requested bigraph exists if
and only if the intersection of $Q'_S$ and the plank $K_S(\alpha
,\beta )$ in ${\bf R}\sp S$ is non-empty.  By Claim \ref{savmetszet}
this intersection is non-empty if and only if $p'(S) \leq \beta $ and
$\alpha \leq b'(S)$.

By applying \eqref{(minb)} to $U=S$ and $U'=X\cup Z$ (where $X\subseteq
S$, \ $Z\subseteq T$), we obtain that $\alpha \leq b'(S)$ is
equivalent to requiring $$ \hbox{ $\alpha \leq [\vert T-Z\vert \vert
X\vert - \widetilde p_T({\cal T}) + \vert {\cal T}\vert \vert X\vert ]
- \widetilde f'(Z) + \widetilde {g'}(S-X)$ }\ $$ whenever $X\subseteq
S, Z\subseteq T,$ and ${\cal T}$ is a subpartition of $Z$.  By letting
$Y=T-Z$ and observing that $ - \widetilde f'(Z) + \widetilde {g'}(S-X)
= \widetilde g_T(T-Y) + \widetilde g_S(S-X)$, we conclude that $\alpha
\leq b'(S)$ is equivalent to \eqref{(galfa)}.

Let $U=S$, \ $Y= T-Z, \ X=S-X'$, \ $U'=X'\cup Y$.  \ Then $V-U'=X\cup
Z, \ U'-S=Y, \ S-U'=X$, and $p_0(U') = - b_0(V-U') = -b_0(X\cup Z)=
-b_0(X\cup (T-Y))$.  Furthermore $$p'(S) = \max \{p_0(U') - \widetilde
{g'}(U'-S) + \widetilde f'(S-U'):  \ U'\subseteq V\}= $$ $$ = \max
\{-b_0(X\cup (T-Y)) + \widetilde f_T(Y) + \widetilde f_S(X):  \
X\subseteq S, Y\subseteq T\}.  $$ Hence $\beta \geq p'(S)$ is
equivalent to $$\beta \geq - [\vert Y\vert \vert X\vert
-\widetilde p_T({\cal T}) +\vert {\cal T}\vert \vert X\vert ] +
\widetilde f_T(Y) + \widetilde f_S(X)$$ $$\text{for every $X\subseteq S,
Y\subseteq T$, $\cal T$ a subpartition of $T-Y$,}$$ and this is
just \eqref{(fbeta)}.  $\bullet$\endproof

\begin{corollary}Provided that there is a simple bigraph covering $p_T$ and
degree-constrained by $(f,g)$, the minimum number of edges of such a
bigraph is \begin{eqnarray}&\max \{ \widetilde f_S(X) + \widetilde f_T(Y) -
\vert X\vert \vert Y\vert + \widetilde p_T({\cal T}) - \vert X\vert
\vert {\cal T}\vert : & \nonumber\\[5pt] & X\subseteq S, \ Y\subseteq T,\ {\cal T} \text{a subpartition of $T-Y$}\}.& \label{(maxelszam)} \end{eqnarray} \end{corollary}

Analogous theorem can be formulated for the maximum number of edges,
as well.

\medskip 

\section{Packing branchings and arborescences}

Let $D=(V,A)$ be a digraph on $n$ nodes.  An {\bf arborescence} is a
directed tree in which one node, its root-node, has no entering arc
and the in-degree of all other nodes is 1. A {\bf branching} $(V,B)$
of $D$ is a directed forest consisting of arborescences.  Its {\bf
root-set} $R(B)$ is the set of nodes of in-degree zero.  By the {\bf
size} of a branching we mean the number of its arcs while the {\bf
root-size} is $\vert R(B)\vert $. Obviously, $\vert B\vert +\vert
R(B)\vert =n.$ In what follows the same term $B$ will be used for a
branching and for its set of arcs.

$D$ is called {\bf rooted $k$-edge-connected} with respect to a
root-node $r_0$ if $\varrho _D(X)\geq k$ for every $\emptyset \subset
X\subseteq V-r_0$.  By Menger, this is equivalent to requiring that
there are $k$ edge-disjoint paths from $r_0$ to $v$ for every node
$v\in V$.

\subsection{Background}

A major open problem in combinatorial optimization is to find a good
characterization for the existence of $k$ disjoint common bases of two
matroids.  This is solved only in special cases.  For example, $\mu
$-element matchings of a bipartite graph form the common bases of two
matroids.  Folkman and Fulkerson \cite{FoFu} proved the following.

\begin{theorem}\label{FoFu} A bigraph $G=(S,T;E)$ includes $k$ disjoint
matchings of size $\mu $ if and only if $$ \hbox{ $k (\mu +\vert
Z\vert - \vert S\cup T\vert ) \leq i_G(Z)$ whenever $Z\subseteq S\cup T$, }\
$$ where $i_G(Z)$ denotes the number of edges induced by $Z$.\end{theorem}

As the branchings of a digraph form the common independent sets of two
matroids, the problems of finding $k$ disjoint spanning arborescences
or $k$ disjoint branchings of size $\mu $ can also be viewed as
special cases of the disjoint common bases problem.  This matroidal
aspect particularly underpins the significance of the following
fundamental result of Edmonds \cite{Edmonds73}.

\begin{theorem}[Edmonds] \label{Edmonds} Let $D=(V,A)$ be a digraph.  \

\noindent {\bf (Weak form)} \ $D$ includes $k$ disjoint spanning
arborescences with a specified root-node $r_0$ if and only if $D$ is
rooted $k$-edge-connected.

\noindent {\bf (Strong form)} \ $D$ includes $k$ disjoint branchings
with specified root-sets $R_1,R_2,\dots ,R_k$ if and only if $\varrho
_D(X)\geq p_R(X)$ for $X\subseteq V$ where $p_R(X)$ denotes the number
of root-sets disjoint from $X$ when $X\not =\emptyset $ and
$p_R(\emptyset )=0$.\end{theorem}

Though Lov\'asz \cite{Lovasz76a} found a short proof relying on
submodular functions and also a great number of variations and
generalizations have been developed (see the book of Schrijver
\cite{Schrijverbook} or a recent survey by Kamiyama
\cite{Kamiyama14}), Edmonds' theorem and the topic of disjoint
branchings remained rather isolated from general frameworks like the
one of submodular flows.  Due to its specific position within
combinatorial optimization, it is particularly important to
investigate extensions and variations.

An early variation of the weak form was proved in \cite{FrankP4}.

\begin{theorem}\label{freeroot} A digraph $D$ has $k$ disjoint spanning
arborescences with unspecified roots {\em (that is, $k$ disjoint
branchings of size $\vert V\vert -1$)} if and only if $$ \hbox{ $\sum
_{i=1}\sp q \varrho _D(V_i)\geq k(q-1)$ for every subpartition
$\{V_1,\dots ,V_q\}$ of $V$.  }\ $$ \end{theorem}

The following extension is due to Cai \cite{Cai83} and Frank
\cite{FrankP4} (see also Theorem 10.1.11 in the book
\cite{Frank-book}).

\begin{theorem}\label{Cai-Frank} Let $f:V\rightarrow {\bf Z}_+$ and
$g:V\rightarrow {\bf Z}_+$ be lower and upper bounds for which $f\leq
g$.  A digraph $D=(V,A)$ includes $k$ disjoint spanning arborescences
so that each node $v$ is the root of at least $f(v)$ and at most
$g(v)$ of these arborescences if and only if $\widetilde f(V)\leq k$,
\begin{equation}\hbox{ $\sum _{i=1}\sp q \varrho _D(V_i)\geq k(q-1) + \widetilde
f(V_0)$ for every partition $\{V_0,V_1,\dots ,V_q\}$ of $V$ }\
\label{(Cai-Frank)} \end{equation} where $q\geq 1$ and only $V_0$ can be empty,
and $$ \hbox{ $\widetilde g(X)\geq k-\varrho _D(X)$ for every subset
$\emptyset \subset X\subseteq V.$}\ $$ \end{theorem}

Note that the condition $\widetilde f(V)\leq k$ can be interpreted as
the inequality in \eqref{(Cai-Frank)} written for $q=0$.  With similar
techniques, the following generalization of Theorem \ref{freeroot} can
also be derived (though, to our best knowledge, it was not explicitly
formulated earlier.)

\begin{theorem}\label{muelu} A digraph $D$ has $k$ disjoint branchings of
size $\mu $ if and only if $$ \hbox{ $\sum _{i=1}\sp q \varrho _D(V_i)
\geq k[q-(n-\mu )] $ for every subpartition $\{V_1,\dots ,V_q\}$ of
$V.$}\ $$ \end{theorem}

\subsection{Packing branchings with prescribed sizes}

The following possible extension emerges naturally for branchings and
matchings, as well.  What is a necessary and sufficient condition for
the existence of $k$ disjoint branchings in a digraph (respectively,
$k$ disjoint matchings in a bigraph) having prescribed sizes $\mu
_1,\mu _2,\dots ,\mu _k$?  A bit surprisingly, the answer in the two
cases is quite different.  For bipartite matchings the problem was
shown to be {\bf NP}-complete even for $k=2$ (\cite{KaMk},
\cite{Palvolgyi}, \cite{Puleo}).  On the other hand, for branchings we
have the following straight generalization of Theorem \ref{muelu}.

\begin{theorem}\label{arbo.basic} Given $k$ positive integers $\mu _1,\mu
_2,\dots ,\mu _k$ $(\mu _j\leq n-1)$, a digraph $D=(V,A)$ on $n$ nodes
has $k$ disjoint branchings $B_1,\dots ,B_k$ of sizes $\vert B_j\vert
=\mu _j$ \ $(j=1,\dots ,k)$ if and only if \begin{equation}\hbox{ $\sum _{i=1}\sp
q \varrho _D(V_i) \geq \sum _{j=1}\sp k [q-(n-\mu _j)]\sp + $ for
every subpartition ${\cal P}= \{V_1,\dots ,V_q\}$ of $V$.  }\
\label{(arbo.main)} \end{equation} \end{theorem}

\proof{Proof.} Throughout we use the notation $m_j:=n-\mu _j$.

Necessity.  The root-set $R_j$ of a branching $B_j$ of size $\mu _j$
has $m_j$ elements.  If $B_j$ has no arc entering $V_i$, then $R_j$
has an element in $V_i$, therefore there are at least $(q-m_j)\sp +$
arcs of $B_j$ entering a member of the subpartition ${\cal
P}=\{V_1,\dots ,V_q\}$, implying that the total number $\sum _{i=1}\sp
q \varrho _D(V_i)$ of arcs entering some members of $\cal P$ is at
least $\sum _{j}\sp k [q-(n-\mu _j)]\sp +$.  (Note that the assumption
$(\mu _j\leq n-1)$ is actually superfluous since \eqref{(arbo.main)},
when applied to $q=1$ and ${\cal P}=\{V\}$, implies that $0=\varrho
_D(V) \geq \sum _{j=1}\sp k [1-(n-\mu _j)]\sp +$ from which each
summand $[1-(n-\mu _j)]\sp +$ must be zero, that is, $1\leq n-\mu
_j$.)

Sufficiency.  Let $S=\{s_1,s_2,\dots ,s_k\}$ be a set of $k$ elements.
We may consider $S$ as the index set of the $k$ branchings to be
found.  Define $m_S:S\rightarrow {\bf Z}_+$ by $m_S(s_j):=m_j$ \
$(j=1,\dots ,k)$.  Let $T:=V$ and define a set-function $p_T$ on $T$
as follows.  \begin{equation}p_T(Y):= \begin{cases}
k-\varrho _D(Y) & \text{if $\emptyset\subset Y\subseteq T$}\\ 0 & \text{if $Y=\emptyset $.} \end{cases} \label{(pTdef)} \end{equation}

Then $p_T$ is intersecting supermodular.  From \eqref{(arbo.main)}, we
have $$\sum _{i=1}\sp q \varrho _D(V_i) \geq \sum _{j=1}\sp k
(q-m_j)\sp + = \sum _{j=1}\sp k \max \{q-m_j,0\}= kq + \sum _{j=1}\sp
k \max \{-m_j,-q\} = kq - \sum _{j=1}\sp k \min \{m_j,q\}$$ from which
$$ \sum _{i=1}\sp q p_T(V_i) = \sum _{i=1}\sp q [k-\varrho _D(V_i)]
\leq \sum _{j=1}\sp k \min \{m_j,q\} = \sum _{s\in S} \min
\{m_S(s),q\}.$$ Therefore \eqref{(basic.felt.free)} holds and Theorem
\ref{basic.base.free} implies that there is a simple bigraph
$G=(S,T;E)$ covering $p_T$ for which $d_G(s)=m_S(s)$ for every $s\in
S$.

For $s_j\in S$ let $R_j$ denote the set of neighbours of $s_j$ in $G$.
Then $\vert R_j\vert =m_j$ for $j=1,\dots ,k$.  Since each non-empty
subset $Y$ of $V$ has at least $p_T(Y)=k-\varrho _D(Y)$ neighbours,
the number of non-neighbours is at most $\varrho _D(Y)$, that is, the
number of sets $R_j$'s disjoint from $Y$ is at most $\varrho _D(Y)$.
The strong form of Edmonds' theorem implies that there are $k$
disjoint branchings $B_1,\dots ,B_k$ with root sets $R_1,\dots ,R_k$,
respectively.  By the definition of $m_j$, we have $\vert B_j\vert = n
- \vert R_j\vert = n - m_j = \mu _j$.  $\bullet$\endproof

\medskip With a similar approach, we can characterize the situation
when not only the sizes of the $k$ disjoint branchings are specified
but the indegree of each node in their union, as well.

\begin{theorem}\label{arbo.mSmT} Let $D=(V,A)$ be a digraph on $n$ nodes,
$m_{\rm in}:V\rightarrow {\bf Z}_+$ an in-degree prescription with
$0\leq m_{\rm in}(v)\leq \varrho _D(v)$ and $m_{\rm in}(v)\leq k$ for
each $v\in V$.  Let $\mu _1,\mu _2,\dots ,\mu _k$ be $k$ positive
integers such that $\mu _1+\cdots +\mu _k=\widetilde m_{\rm in}(V)$.
There is a subgraph $(V,F)$ of $D$ which is the union of $k$ disjoint
branchings $B_1,\dots ,B_k$ of sizes $\vert B_j\vert =\mu _j$ \
$(j=1,\dots ,k)$ and for which $$ \hbox{ $\varrho _F(v)=m_{\rm in}(v)$
for each $v\in V$}\ $$ if and only if \begin{equation}\widetilde m_{\rm in}(Y) +
\sum _{i=1}\sp q \varrho _D(V_i) \geq \sum _{j=1}\sp k [q+\vert Y\vert
-(n-\mu _j)]\sp + \label{(arbo.mSmT)} \end{equation} for every subset
$Y\subseteq V$ and every subpartition $\{V_1,\dots ,V_q\}$ of $V-Y$.
\end{theorem}

\proof{Proof.} Necessity.  Suppose that the requested $k$ branchings
$B_1,\dots ,B_k$ exist and let $F=B_1\cup \cdots \cup B_k$.  Let
$Y\subseteq V$ and ${\cal P}=\{V_1,\dots ,V_q\}$ be a subpartition of
$V-Y$.  As before, $m_j=n-\mu _j$ is the cardinality of the root-set
$R_j$ of $B_j$.  Therefore the number of non-root nodes in $Y$
($=\vert Y-R_j\vert )$ plus the number of $V_i$'s disjoint from $R_j$
is at least $\vert Y\vert +q - m_j$, and hence the number of arcs of
$B_j$ entering a node of $Y$ plus the number of arcs of $B_j$ entering
a member of $\cal P$ is at least $(\vert Y\vert +q - m_j)\sp +$.
Hence $$ \widetilde m_{\rm in}(Y) + \sum _{i=1}\sp q \varrho _D(V_i)
\geq \widetilde m_{\rm in}(Y) + \sum _{i=1}\sp q \varrho _F(V_i)=$$ $$
\sum _{j=1}\sp k\ [\sum _{v\in Y} \varrho _{B_j}(v) + \sum _{i=1}\sp q
\varrho _{B_j} (V_i)] \geq \sum _{j=1}\sp k (\vert Y\vert +q - m_j)\sp
+,$$ and \eqref{(arbo.mSmT)} follows.

Sufficiency.  Let $S,T$, and $m_S$ be the same as in the preceding
proof.  Define a set-function $p_T$ on $T$ as follows.  \begin{equation}p_T(Y):=
\begin{cases}
k-\varrho _D(Y) & \text{if $Y\subseteq T, \ \vert Y\vert \geq 2$}\\ k-m_{\rm in}(v) & \text{if $Y=\{v\},\ v\in V$}\\ 0 & \text{if $Y=\emptyset
$.}\end{cases} \label{(pTdefY)} \end{equation} The hypothesis $m_{\rm in}(v)\leq
\varrho _D(v)$ implies that $k-m_{\rm in}(v)\geq k-\varrho _D(v)$ and
hence $p_T$ is intersecting supermodular.  Let ${\cal T}=\{V_1,\dots
,V_q,V_{q+1},\dots ,V_{q'}\}$ be a subpartition of $V$ so that the
first $q$ members are of cardinalities at least two while the
subsequent members are singletons.  Let ${\cal P}=\{V_1,\dots ,V_q\}$
and let $Y$ denote the union of $V_{q+1},\dots ,V_{q'}$ (that is,
$\vert Y\vert = q'-q$).

By \eqref{(arbo.mSmT)}, we have $$ \widetilde m_{\rm in}(Y) + \sum
_{i=1}\sp q \varrho _D(V_i) \geq \sum _{j=1}\sp k [q'-(n-\mu _j)]\sp +
= \sum _{j=1}\sp k \max \{q'-m_j,0\} = $$ $$ kq' + \sum _{j=1}\sp k
\max\{-m_j,-q\} = kq' - \sum _{j=1}\sp k \min\{m_j,q'\}$$ from which
$$ \sum _{i=1}\sp {q'} p_T(V_i) = \sum _{v\in Y} [k-m_{\rm in}(v)] +
\sum _{i=1}\sp q [k-\varrho _D(V_i)] = k(\vert Y\vert +q) - \widetilde
m_{\rm in}(Y) - \sum _{i=1}\sp q \varrho _D(V_i) =$$ $$ kq' -
\widetilde m_{\rm in}(Y) - \sum _{i=1}\sp q \varrho _D(V_i) \leq \sum
_{j=1}\sp k \min\{m_j,q'\} = \sum _{j=1}\sp k \min\{m_S(s_j),q'\}.  $$
Therefore \eqref{(basic.felt.free)} holds with $q'$ in place of $q$ and
with $V_i$ in place of $T_i$, and Theorem \ref{basic.base.free}
implies that there is a simple bigraph $G=(S,T;E)$ covering $p_T$ for
which $d_G(s)=m_S(s)$ for every $s\in S$.  Since $G$ covers $p_T$, it
follows that $d_G(v)\geq p_T(v) = k - m_{\rm in}(v)$.  Hence $$ \sum
_{v\in V} [k - m_{\rm in}(v)] \leq \sum _{v\in V}d_G(v) = \sum _{s\in
S}d_G(s) = \sum _{s\in S}m_S(s) = \sum _{j=1}\sp k (n- \mu _j).$$ But
here we must have equality since we assumed that $\mu _1+\cdots +\mu
_k=\widetilde m_{\rm in}(V)$.  This implies that $d_G(v) = k - m_{\rm
in}(v)$ for each $v\in V$.

For $s_j\in S$ let $R_j$ denote the set of neighbours of $s_j$ in $G$.
Then $\vert R_j\vert =m_j$ for $j=1,\dots ,k$.  Since each non-empty
subset $Y$ of $V$ has at least $p_T(Y)=k-\varrho _D(Y)$ neighbours,
the number of non-neighbours is at most $\varrho _D(Y)$, that is, the
number of sets $R_j$'s disjoint form $Y$ is at most $\varrho _D(Y)$.
The strong form of Edmonds' theorem implies that there are $k$
disjoint branchings $B_1,\dots ,B_k$ with root sets $R_1,\dots ,R_k$,
respectively.  By the definition of $m_j$, we have $\vert B_j\vert = n
- \vert R_j\vert = n - m_j = \mu _j$.

Let $F:=B_1\cup \cdots \cup B_k$.  As $d_G(v)$ is the number of
$R_j$'s containing $v$, the indegree $\varrho _F(v)$ is $k-d_G(v) =
m_{\rm in}(v)$, as required.  $\bullet$\endproof

\medskip Note that the indegree $\varrho _F(v)$ in the union $F$ of
$k$ disjoint branchings is exactly $k$ minus the number of root-sets
not containing $v$.  Therefore Theorem \ref{arbo.mSmT} could be
described in an equivalent form when, instead of the indegree of each
node $v$ in the union of $k$ branchings with specified sizes, the
number of root-sets containing $v$ is prescribed.

\subsection{Packing branchings with bounds on sizes, on total
indegrees, and on total size}

Suppose now that, instead of exact prescription $\mu _j$ for the size
of the branchings $B_j$, we are given a lower bound $\varphi _j$ and
an upper bound $\gamma _j$ with $0\leq \varphi _j\leq \gamma _j\leq
n-1$ \ $(j=1,\dots ,k)$.  Furthermore, instead of the exact
prescription $m_{\rm in}(v)$ for the indegree $\varrho _F(v)$ \ ($v\in
V$), where $F$ denotes the union of the $k$ branchings, we are given a
lower bound $f_{\rm in}(v)$ and an upper bound $g_{\rm in}(v)$ for
which $0\leq f_{\rm in}(v) \leq g_{\rm in}(v) \leq k$.  Moreover, we
impose a lower bound $\alpha _{\rm u}$ \ and an upper bound $\beta
_{\rm u}$ for the cardinality of the union of the $k$ branchings.

The proof of Theorem \ref{arbo.basic} relied on a one-to-one
correspondence between simple bigraphs $G=(S,T;E)$ covering the
function $p_T$ defined in \eqref{(pTdef)} \ (where $T=V$ and
$S=\{s_1,\dots ,s_k\}$ is a $k$-element index set of the $k$
branchings to be found) and the families ${\cal R}=\{R_1,\dots ,R_k\}$
of $k$ root-sets satisfying the necessary condition in the strong form
of Edmonds' theorem (which required that $\varrho _D(Y)$ is at least
the number of $R_i$'s disjoint from $Y$ for each non-empty $Y\subseteq
V$).  Let $B_1,\dots ,B_k$ denote the $k$ disjoint branchings ensured
by Edmonds' theorem for which $R(B_j)=R_j$, and let $F=B_1\cup \cdots
\cup B_k$.

In this correspondence, the degree of a node $s_j\in S$ is the
cardinality of $R_j$, that is, $$d_G(s_j)=\vert R_j\vert = n - \vert
B_j\vert .$$ Furthermore, the degree of a node $v\in V=T$ is the
number of root-sets $R_j$'s containing $v$, that is, $$d_G(v)=
k-\varrho _F(v).$$ Finally, for the total number of edges of $G$, we
have $$\vert E\vert =\sum _{j=1}\sp k d_G(s_j) = \sum _{j=1}\sp k
\vert R_j\vert = \sum _{j=1}\sp k (n-\vert B_j\vert )=nk -\vert F\vert
.$$ Define $$ \hbox{ $f_S(s_j):=n-\gamma _j$ and $g_S(s_j):=n-\varphi
_j$ \ for $s_j\in S$, }\ $$ $$ \hbox{ $f_T(v):=k - g_{\rm in}(v)$ and
$g_T(v):=k - f_{\rm in}(v)$ \ for $v\in T=V$, }\ $$ $$ \hbox{ $\alpha
:= kn - \beta _{\rm u}$ and $\beta := kn - \alpha _{\rm u}$.  }\ $$

By this vocabulary, $\varphi _j\leq \vert B_j\vert \leq \gamma _j$ if
and only if $f_S(s_j)\leq d_G(s_j) \leq g_S(s_j)$ \ $(s_j\in S).$
Furthermore, $f_{\rm in}(v) \leq \varrho _F(v)\leq g_{\rm in}(v)$ if
and only if $f_T(v)\leq d_G(v) \leq g_T(v)$ \ $(v\in T)$.  Finally,
$\alpha _{\rm u} \leq \vert B_1\cup \cdots \cup B_k\vert \leq \beta
_{\rm u}$ if and only if $\alpha \leq \vert E\vert \leq \beta $. By
aggregating Theorems \ref{fg.main} and \ref{fg.ab}, we obtain the
following.

\begin{theorem}In a digraph $D=(V,A)$ on $n$ nodes, there are $k$ disjoint
branchings $B_1,\dots ,B_k$ for which $\varphi _j\leq \vert B_j\vert
\leq \gamma _j$ \ $(j=1,\dots ,k)$, for which $f_{\rm in}(v)\leq
\varrho _F(v) \leq g_{\rm in}(v)$ \ $(v\in V)$, and for which $\alpha
_{\rm u}\leq \vert F\vert \leq \beta _{\rm u}$, where $F=B_1\cup
\cdots \cup B_k$, if and only if the conditions \eqref{(ftgs)},
\eqref{(fsgt)}, \eqref{(galfa)}, and \eqref{(fbeta)} hold for the
choice of $f_T,g_T,f_S,g_S, \alpha , \beta $ defined above.  $\bullet
$ \end{theorem}

\medskip 

\section{Maximum term rank problems} \label{termrank}

\subsection{Degree-specified max term rank}

The members of ${\cal G}(m_S,m_T)$ (that is, simple bigraphs fitting
the degree-specification $(m_S,m_T)$) can be identified with
$(0,1)$-matrices of size $\vert S\vert \vert T\vert $ with row sum
vector $m_S$ and column sum vector $m_T$.  Let ${\cal M}(m_S,m_T)$
denote the set of these matrices.  Ryser \cite{Ryser58} defined the
{\bf term rank} of a $(0,1)$-matrix $M$ by the maximum number of
independent $1$'s which is the matching number of the bipartite graph
corresponding to $M$.  Ryser developed a formula for the maximum term
rank of matrices in ${\cal M}(m_S,m_T)$.  The maximum term rank
problem is equivalent to finding a bipartite graph $G$ in ${\cal
G}(m_S,m_T)$ whose matching number $\nu (G)$ is as large as possible.
Although we use graph terminology, the original name \lq term rank\rq
\ for the problem will be kept throughout.  In graphical terms,
Ryser's theorem is equivalent to the following.

\begin{theorem}[Ryser] \label{term.maxterm} Let $\ell\leq \vert T\vert $ be
an integer.  Suppose that ${\cal G}(m_S,m_T)$ is non-empty, that is,
Condition \eqref{(ujterm.GR)} holds.  Then ${\cal G}(m_S,m_T)$ has a
member $G$ with matching number $\nu (G)\geq \ell$ if and only if \begin{equation}
\widetilde m_S(X) + \widetilde m_T(Y) -\vert X\vert \vert Y\vert +
(\ell - \vert X\vert -\vert Y\vert ) \leq \gamma \ \hbox{whenever}\
X\subseteq S, \ Y\subseteq T. \label{(term.maxterm)} \end{equation} Moreover,
\eqref{(term.maxterm)} holds if the inequality in it is required only
when $X$ consists of the $i$ largest values of $m_S$ and $Y$ consists
of the $j$ largest values of $m_T$ $(i=0,1,\dots ,\vert S\vert , \
j=0,1,\dots ,\vert T\vert )$.  \end{theorem}

Observe that the conditions \eqref{(term.maxterm)} and
\eqref{(ujterm.GR)} in Theorem \ref{term.maxterm} can be united as
follows.  \begin{equation}\widetilde m_S(X) + \widetilde m_T(Y) -\vert X\vert
\vert Y\vert + (\ell - \vert X\vert -\vert Y\vert )\sp + \leq \gamma \
\hbox{ whenever}\ X\subseteq S, \ Y\subseteq T
\label{(term.maxterm.uni)}, \end{equation} that is, assuming this inequality, we
do not need to impose explicitly the non-emptiness of ${\cal
G}(m_S,m_T)$.

Note that the strengthening formulated in the second part of the
theorem is nothing but a straightforward observation.  Beyond the
aesthetic joy, a practical advantage is that such simplified condition
can easily be checked in polynomial time since there are only a few
($(\vert S\vert +1)(\vert T\vert +1)$) inequalities to be checked.
This will be crucial in the algorithm described below for the
degree-constrained max term rank problem.  Note that the original
proof of Ryser gives rise to a polynomial time algorithm to compute
the matrix itself.  Subsequently, Brualdi and Ross \cite{Brualdi-Ross}
described a simpler proof and which gives rise to a simple algorithm.

We also remark that there is a characterization given by Haber
\cite{Haber60} for the minimum term rank of the graphs in ${\cal
G}(m_S,m_T)$ but we deal only with the maximum term rank problem.

\subsubsection{Relation to network flows}

As the bipartite matching problem and the more general
degree-prescribed subgraph problem can be treated with network flow
technique, one may be wondering if Ryser's theorem could also be
derived via network flows.  Ford and Fulkerson, for example, remarked
in their classic book (\cite{Ford-Fulkerson}, p. 89) that:

$$ \hbox{ {\it \lq Neither term rank problem appears amenable to flow
approach.\rq} }\ $$

Such a link could help solving the weighted and the subgraph version
of the max term rank problem.  But recently it turned out that the
failure of the attempt of Ford and Fulkerson was not just by chance.
It was proved (\cite{KaMk}, \cite{Palvolgyi}, \cite{Puleo}) that the
problem of deciding whether an initial bigraph $G_0$ has a perfectly
matchable degree-specified subgraph is {\bf NP}-complete.  Therefore
both the weighted and the subgraph versions of the max term rank
problem is {\bf NP}-complete, showing that even the theory of
submodular flows cannot help.  The first goal of this section is to
show that Ryser's theorem immediately follows from the general result
on covering a supermodular function by simple bigraphs developed in
Section \ref{ujbasic}.

Unfortunately not only weighted, but quite natural unweighted
extensions also turned out to be {\bf NP}-complete.  For example,
finding a member of ${\cal G}(m_S,m_T)$ in which there is a subgraph
with specified degrees is equivalent to finding two disjoint simple
bipartite graphs on the same node-set, and this latter problem was
shown to be {\bf NP}-complete by D\"urr, Guinez, and Matamala
\cite{DGM} (see Proposition \ref{DGMa}).

In the light of the {\bf NP}-complete problems in the close
neighbourhood, it is pleasing to realize that there are nicely
tractable extensions of the max term rank problem.  In the present
section, we shall extend Ryser's theorem to the case when the bigraph
with high matching number is degree-constrained and edge-number
constrained, not just degree-specified.

In paper \cite{Berczi-Frank16b}, we shall develop an augmentation and
a matroidal generalization.  In the first one, a given initial bigraph
is to be augmented to get a simple degree-specified bigraph with
matching number at least $\ell.$ In matrix terms, this means that some
of the entries of the $(0,1)$-matrix are specified to be 1. This is in
sharp contrast with the {\bf NP}-completeness of that version when
some entries of the matrix are specified to be 0. In the matroidal
extension of Ryser's theorem, there are matroids on $S$ and on $T$ and
we want to find a degree-specified simple bigraph including a matching
that covers bases in both matroids.

\subsubsection{Proof of Ryser's theorem}

\proof{Proof.} Necessity.  Let $G$ be a bipartite graph with the requested
properties.  Since $G$ is simple, it has at least $\widetilde m_S(X) +
\widetilde m_T(Y) -\vert X\vert \vert Y\vert $ edges having at least
one end-node in $X\cup Y$.  Moreover, since $G$ has a matching of
$\ell$ edges, there are at least $\ell -\vert X\cup Y\vert $ edges
connecting $S-X$ and $T-Y$.  Therefore the total number $\gamma $ of
edges is at least $\widetilde m_S(X) + \widetilde m_T(Y) -\vert X\vert
\vert Y\vert + \ell -\vert X\cup Y\vert $, that is,
\eqref{(term.maxterm)} is indeed necessary.

Sufficiency.  We need the following deficiency form of Hall's theorem.

\begin{lemma}[Hall and Ore] \label{Hall-Ore} Let $G=(S,T;E)$ be a bipartite
graph and $\ell\leq \vert T\vert $ an integer.  The matching number
$\nu (G)$ is at least $\ell $ {\em (that is, there is a matching of
$\ell $ edges)} if and only if \begin{equation}\vert \Gamma (Y)\vert \geq \ell
-\vert T-Y\vert \ \hbox{ holds for every }\ Y\subseteq T.
\label{(termrank.Ore)} \end{equation} \end{lemma}

Define a set-function $p_T$ on $T$ by \begin{equation}p_T(Y):= \begin{cases}
\ell -(\vert T-Y\vert ) & \text{if $\emptyset \subset Y\subseteq T$} \\ 0 & \text{if
$Y=\emptyset$.}\end{cases} \label{(termrank.pTdef)} \end{equation} Then $p_T$ is
fully supermodular and monotone non-decreasing.  Since ${\cal
G}(m_S,m_T)$ is assumed to be non-empty, \eqref{(ujbasic.nemures1)}
holds.  By \eqref{(term.maxterm)}, we have $$\widetilde m_S(X) +
\widetilde m_T(Y) -\vert X\vert \vert Y\vert + p_T(T-Y) - \vert X\vert
= \widetilde m_S(X) + \widetilde m_T(Y) -\vert X\vert \vert Y\vert +
\ell -\vert Y\vert - \vert X\vert \leq \gamma ,$$ that is,
\eqref{(ujbasic.nemures2)} also holds.  By Corollary \ref{msmtfullp}
(in Section \ref{ptst}), there is a simple bipartite graph $G=(S,T;E)$
covering $p_T$ and fitting $(m_S,m_T)$.  But such a graph has a
matching of size $\ell$ by Lemma \ref{Hall-Ore}, and we are done.
$\bullet$\endproof

\subsection{Degree and edge-number constrained max term rank}

Our goal is to extend Ryser's theorem for the case when upper or lower
bounds are given for the degrees rather than exact prescriptions.
Bounds for the total number of edges can also be incorporated.  Let
$f_V=(f_S,f_T)$ and $g_V=(g_S,g_T)$ be lower and upper bound functions
with $0\leq f_V\leq g_V$.  As we are interested in simple bigraphs, we
may suppose that $g_S(s)\leq \vert T\vert $ for every $s\in S$ and
$g_T(t)\leq \vert S\vert $ for every $t\in T$.

Ryser's theorem was derived above by applying Corollary
\ref{msmtfullp} to the set-function $p_T$ defined in
\eqref{(termrank.pTdef)}.  By applying Corollary \ref{bigraph-linking}
to the same $p_T$, we obtain the following extension.

\begin{theorem}\label{termrank.link1} Let $\ell\leq \vert T\vert $ be an
integer, $f_V=(f_S,f_T)$ and $g_V=(g_S,g_T)$ bounds with $f_V\leq
g_V$.  \medskip

\noindent {\bf (A)} \ There is a simple bigraph $G'=(S,T;E')$ with
matching number $\nu (G')\geq \ell$ and degree-constraints $(f_T,g_S)$
if and only if \begin{equation}\widetilde f_T(Y) - \vert X\vert \vert Y\vert +
(\ell - \vert X\vert -\vert Y\vert )\sp + \leq \widetilde g_S(S-X) \
\hbox{ whenever }\ X\subseteq S, \ Y\subseteq T.
\label{(termrank.ftgs)} \end{equation} Moreover, \eqref{(termrank.ftgs)} holds
if the inequality in it is required only when $X$ consists of the $i$
largest values of $g_S$ and $Y$ consists of the $j$ largest values of
$f_T$ $(i=0,1,\dots ,\vert S\vert , \ j=0,1,\dots ,\vert T\vert )$.
\medskip

\noindent {\bf (B)} \ There is a simple bigraph $G''=(S,T;E'')$ with
matching number $\nu (G'')\geq \ell$ and degree-constraints
$(f_S,g_T)$ if and only if \begin{equation}\widetilde f_S(X) - \vert X\vert \vert
Y\vert + (\ell - \vert X\vert -\vert Y\vert )\sp + \leq \widetilde
g_T(T-Y) \ \hbox{ whenever}\ X\subseteq S, \ Y\subseteq T.
\label{(termrank.fsgt)} \end{equation} Moreover, \eqref{(termrank.fsgt)} holds
if the inequality in it is required only when $X$ consists of the $i$
largest values of $f_S$ and $Y$ consists of the $j$ largest values of
$g_T$ $(i=0,1,\dots ,\vert S\vert , \ j=0,1,\dots ,\vert T\vert )$.
\medskip

\noindent {\bf (AB)} \ There is a simple bigraph $G=(S,T;E)$ with
matching number $\nu (G)\geq \ell$ and degree-constraints $(f_V,g_V)$
if and only if both $G'$ and $G''$ exist (that is, both
\eqref{(termrank.ftgs)} and \eqref{(termrank.fsgt)} hold).  $\bullet
$ \end{theorem}

By applying Theorem \ref{fg.ab} to the same $p_T$ defined in
\eqref{(termrank.pTdef)}, we obtain the following extension.

\begin{theorem}\label{termrank.link2} Suppose that there is a simple bigraph
with matching number at least $\ell$ which is degree-constrained by
$(f_V,g_V)$ (that is, conditions \eqref{(termrank.ftgs)} and
\eqref{(termrank.fsgt)} hold).  There is simple bigraph $G=(S,T;E)$
with matching number at least $\ell$ which is degree-constrained by
$(f_V,g_V)$:

\medskip \noindent {\bf (A)} \ for which $\alpha \leq \vert E\vert $
if and only if \begin{equation}\hbox{ $\widetilde g_S(S-X) + \widetilde g_T(T-Y) +
\vert X\vert \vert Y\vert - (\ell -\vert X\vert -\vert Y\vert )\sp +
\geq \alpha $ whenever $X\subseteq S, \ Y\subseteq T, $}\
\label{(termrank.galfa)} \end{equation}

Moreover, \eqref{(termrank.galfa)} holds if the inequality in it is
required only when $X$ consists of the $i$ largest values of $g_S$ and
$Y$ consists of the $j$ largest values of $g_T$ $(i=0,1,\dots ,\vert
S\vert , \ j=0,1,\dots ,\vert T\vert )$.

\medskip

\noindent {\bf (B)} \ $\vert E\vert \leq \beta $ if and only if \begin{equation}
\hbox{ $ \widetilde f_S(X) + \widetilde f_T(Y) - \vert X\vert \vert
Y\vert + (\ell -\vert X\vert -\vert Y\vert )\sp + \leq \beta $
whenever $X\subseteq S, \ Y\subseteq T,$ }\ \label{(termrank.fbeta)}
\end{equation}

Moreover, \eqref{(termrank.fbeta)} holds if the inequality in it is
required only when $X$ consists of the $i$ largest values of $f_S$ and
$Y$ consists of the $j$ largest values of $f_T$ $(i=0,1,\dots ,\vert
S\vert , \ j=0,1,\dots ,\vert T\vert )$.  \medskip

\noindent {\bf (AB)} \ $\alpha \leq \vert E\vert \leq \beta $ if and
only if both \eqref{(termrank.galfa)} and \eqref{(termrank.fbeta)}
hold. \end{theorem}

\subsubsection{Algorithmic aspects}

As already indicated above, the original proof of Ryser is
algorithmic.  Using this (or a simpler algorithm by Brualdi and Ross
\cite{Brualdi-Ross}) as a subroutine, we describe an algorithm to find
a degree-constrained bigraph with matching number at least $\ell $. A
specific feature of the algorithm is that it makes use of Theorem
\ref{termrank.link1} (and does not re-prove it).  Another basic
constituent is the observation that conditions \eqref{(termrank.ftgs)}
and \eqref{(termrank.fsgt)} can easily be checked in polynomial time,
as stated in the theorem, since it suffices to check the inequalities
in question only for $(\vert S\vert +1)(\vert T\vert +1)$ cases.  The
algorithm starts by checking \eqref{(termrank.ftgs)} and
\eqref{(termrank.fsgt)}, and terminates if anyone of them fails to
hold.  Suppose now that both conditions do hold.

Assume that there is a {\bf loose} node $v$ meaning that
$f_V(v)<g_V(v)$.  We can check in polynomial time whether $f_V(v)$ can
be increased by 1 without destroying \eqref{(termrank.ftgs)} and
\eqref{(termrank.fsgt)}, and if it can, increase $f_V(v)$ by 1. By
repeating this operation as long as possible, we arrive at a situation
where $f_V(v)$ cannot be increased any more at any loose node.

By Theorem \ref{termrank.link1}, there is a simple bigraph $G$ with
$\nu (G)\geq \ell$ and degree-constrained by $(f_V,g_V)$.  Then
$d_G(v) = f_V(v)$ clearly holds for a node with $f_V(v)=g_V(v)$, but
$d_G(v) = f_V(v)$ holds for a loose node $v$, as well, since if we had
$f_V(v)<d_G(v)$, then $f_V(v)$ could be increased without destroying
the conditions.  We can conclude that $m_V:=f_V$ and $\gamma :=f_S(S)$
satisfy \eqref{(term.maxterm)} and therefore Ryser's algorithm (or the
simpler algorithm by Brualdi and Ross) can be applied to construct the
requested $G$.

\medskip \medskip

The same approach works in the case when, in addition to the
degree-constraints $(f_V,g_V)$, there is a lower bound $\alpha $ and
an upper bound $\beta $ for the number of edges.

First, we can check in polynomial time if each of conditions
\eqref{(termrank.ftgs)}, \eqref{(termrank.fsgt)},
\eqref{(termrank.galfa)}, and \eqref{(termrank.fbeta)} holds.  If any of
them is violated, the algorithm terminates.  Suppose that these
conditions hold.  We can also check in polynomial time if there is a
loose node $v$ for which $f_V(v)$ can be increased by 1 without
violating any of these conditions, and we make these liftings of $f_T$
as long as possible.  Therefore the final $f_V$ and $g_V$ continue to
meet the four conditions.  By Theorem \ref{termrank.link1}, there is a
bigraph $G$ satisfying the requirements.

By Theorem \ref{termrank.link2}, there is a simple bigraph $G=(S,T;E)$
with $\nu (G)\geq \ell$ and $\alpha \leq \vert E\vert \leq \beta $
which is degree-constrained by $(f_V,g_V)$.  Then $d_G(v) = f_V(v)$
clearly holds for a node with $f_V(v)=g_V(v)$, but $d_G(v) = f_V(v)$
holds for a loose node $v$, as well, since if we had $f_V(v)<d_G(v)$,
then $f_V(v)$ could be increased without destroying the conditions.
We can conclude that $m_V:=f_V$ and $\gamma :=f_S(S)$ satisfy
\eqref{(term.maxterm)}.

\medskip With a little care, it can be shown that the complexity of
the algorithm above to construct the degree-specification $m_V$
satisfying \eqref{(term.maxterm)} for which $f_V\leq m_V\leq g_V$ and
$\alpha \leq \widetilde m_S(S)\leq \beta $ is $O(n\sp 2\log n)$.

\medskip \medskip

\subsection{Further matching-type requirements}

A special case of the max term rank problem characterizes
degree-specifications which can be realized by a perfectly matchable
bipartite graph.  Brualdi \cite{Brualdi80} characterized
degree-specifications which can be realized by elementary bipartite
graphs.  (A simple bigraph is {\bf elementary} if it is connected and
each of its edges belongs to a perfect matching.)  His result is
extended in \cite{Berczi-Frank16c} to so-called $k$-elementary
bigraphs.

\begin{figure}[t]
\centering
\includegraphics[width=0.6\textwidth]{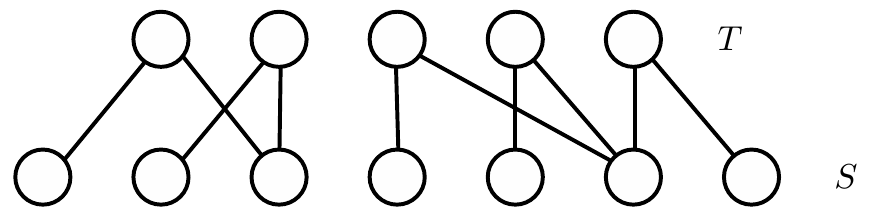}
\caption{A $T_2$-forest}
\label{fig:t2forest}
\end{figure}

In this section, we describe yet another model for degree-specified
bigraphs.  By a {\bf $T_2$-forest} we mean a bigraph $(S,T;F)$ which
is a forest with $d_F(t)=2$ for every $t\in T$ (see Figure~\ref{fig:t2forest}).  Lov\'asz originally
developed Theorem \ref{intro.Lovasz} to characterize bigraphs
$G_0=(S,T;E_0)$ including a $T_2$-forest.

\begin{theorem}\label{Lovasz1m} In a bigraph $G=(S,T;E)$, there exists a
$T_2$-forest if and only if \begin{equation}\vert \Gamma _G(Y)\vert \geq \vert
Y\vert + 1 \ \hbox{whenever }\ \emptyset \not =Y\subseteq T.
\label{(lovasz1b)} \end{equation} \end{theorem}

Lov\'asz used this result to prove a conjecture of Erd{\H o}s on
2-colourability of hypergraphs with the strong Hall inequality.  Here
we show another utilization.

\begin{theorem}\label{wooded} Let $S$ and $T$ be disjoint sets with $\vert
S\vert \geq \vert T\vert +1$ and let $V=S\cup T$.  Let $m_V=(m_S,m_T)$
be a degree-specification for which $\widetilde m_S(S)=\widetilde
m_T(T)=\gamma $ and $m_T(t)\geq 2$ for every $t\in T$.  There exists a
simple bigraph $G=(S,T;E)$ fitting $m_V$ and including a $T_2$-forest
if and only if \eqref{(ujterm.GR)} holds and

\begin{equation}\widetilde m_S(X) + \widetilde m_T(Y) -\vert X\vert \vert Y\vert
-\vert X\vert -\vert Y\vert + \vert T\vert +1 \leq \gamma \ \hbox{
whenever}\ \emptyset \not =X\subseteq S, \ Y\subseteq T.
\label{(T2fa)} \end{equation} \end{theorem}

\proof{Proof.} Necessity.  \ Theorem \ref{ujterm.GR} stated that
\eqref{(ujterm.GR)} was the necessary and sufficient condition for the
realizability of $(m_S,m_T)$.

Suppose that there is a simple bigraph $G=(S,T;E)$ realizing $m_V$ and
including a $T_2$-forest $F$.  The graph has at least $\widetilde
m_S(X) + \widetilde m_T(Y) -\vert X\vert \vert Y\vert $ edges with at
least one end in $X\cup Y$.  Forest $F$ has exactly $2\vert T-Y\vert $
edges ending in $T-Y$.  Among these edges, at most $\vert X\vert +
\vert T-Y\vert -1$ are induced by $X\cup (T-Y)$ since $F$ is a forest
(and $X$ is non-empty).  Therefore $F$ has at least $$2\vert T-Y\vert
- ( \vert X\vert + \vert T-Y\vert -1) = \vert T\vert -\vert X\vert
-\vert Y\vert + 1$$ edges connecting $T-Y$ and $T-X$.  By combining
these observations, we conclude that the left hand side of the
inequality in \eqref{(T2fa)} is indeed a lower bound for the number
$\gamma $ of edges of $G$.

Sufficiency.  \ Define a set-function $p_T$ on $T$ by \begin{equation}p_T(Y):=
\begin{cases} 
\vert Y\vert +1 & \text{if $\emptyset \subset Y\subseteq T$}\\
0 & \text{if $Y=\emptyset.$} \end{cases} \label{(termrank.pTdef1)} \end{equation} Then $p_T$
is intersecting supermodular and monotone non-decreasing.  If there is
a simple bigraph $G$ covering $p_T$ and fitting $m_V$, then Theorem
\ref{Lovasz1m} implies that $G$ has a $T_2$-forest and we are done.
If no such $G$ exists, then Theorem \ref{ujterm.main.spec} implies
that there are subsets $X\subseteq S$, $Y\subseteq T$ and a
subpartition ${\cal T}=\{T_1,\dots ,T_q\}$ of $T-Y$ for which
$$\widetilde m_S(X) + \widetilde m_T(Y) -\vert X\vert \vert Y\vert +
\widetilde p_T({\cal T}) - \vert {\cal T}\vert \vert X\vert >\gamma
,$$ that is, \begin{equation}\widetilde m_S(X) + \widetilde m_T(Y) -\vert X\vert
\vert Y\vert + \sum _{i=1}\sp q (\vert T_i\vert +1) -q\vert X\vert
>\gamma \label{(baj)} \end{equation} We cannot have $q=0$, that is, $\cal T$
cannot be empty because of \eqref{(ujterm.GR)}.

We cannot have $X=\emptyset $, for otherwise $$ \widetilde m_T(Y) +
\vert T-Y\vert + \vert T-Y\vert \geq \widetilde m_T(Y) + \sum
_{i=1}\sp q (\vert T_i\vert +1) = $$ $$\widetilde m_S(X) + \widetilde
m_T(Y) -\vert X\vert \vert Y\vert + \sum _{i=1}\sp q (\vert T_i\vert
+1) -q\vert X\vert >\gamma ,$$ from which $2\vert T-Y\vert > \gamma -
\widetilde m(T)= \widetilde m_T(T-Y),$ contradicting the hypothesis
that $m_T(t)\geq 2$ for each $t\in T$.  Therefore $X\not =\emptyset $.
Since $q\geq 1$ and $\vert X\vert \geq 1$, we have $q(\vert X\vert
-1)\geq \vert X\vert -1$ from which $$ -\vert X\vert -\vert Y\vert +
\vert T\vert +1 = \vert T-Y\vert - (\vert X\vert -1) \geq \vert
T-Y\vert -q(\vert X\vert -1) = $$ $$\vert T-Y\vert + q - q\vert X\vert \geq
\sum _{i=1}\sp q (\vert T_i\vert +1) -q\vert X\vert .$$ This and
\eqref{(baj)} imply $$ \widetilde m_S(X) + \widetilde m_T(Y) -\vert
X\vert \vert Y\vert -\vert X\vert -\vert Y\vert + \vert T\vert +1 \geq
\widetilde m_S(X) + \widetilde m_T(Y) -\vert X\vert \vert Y\vert +
\sum _{i=1}\sp q (\vert T_i\vert +1) -q\vert X\vert >\gamma ,$$
contradicting \eqref{(T2fa)}.  $\bullet$\endproof

\medskip

Actually, Theorem \ref{intro.Lovasz} implies the following more
general form of Theorem \ref{Lovasz1m} in which the forest has a
specified degree (not necessarily identically 2) at each node in $T$.

\begin{theorem}\label{Lovasz1mb} Let $m_{\rm for}:T\rightarrow {\bf Z}_+$ be
a degree specification.  In a bigraph $G=(S,T;E)$, there exists a
forest $F$ with $d_F(t)= m_{\rm for}(t)$ \ $(t\in T)$ if and only if
\begin{equation}\vert \Gamma _G(Y)\vert \geq \widetilde m_{\rm for}(Y)-\vert
Y\vert + 1 \ \hbox{ whenever}\ \emptyset \not =Y\subseteq T.
\label{(lovasz1bx)} \end{equation} \end{theorem} \medskip

Consequently, Theorem \ref{wooded} can also be generalized in such a
way that the simple bigraph should fit a degree specification $m_V$
and should include a forest with specified degrees in the nodes in
$T$.

\begin{theorem}\label{woodedx} Let $S$ and $T$ be disjoint sets and
$V:=S\cup T$.  Let $m_V=(m_S,m_T)$ be a degree-specification for which
$\widetilde m_S(S)=\widetilde m_T(T)=\gamma $, and let $m_{\rm
for}:T\rightarrow {\bf Z}_+$ be a degree specification on $T$ for
which $m_{\rm for} \leq m_T.$ There exists a simple bigraph
$G=(S,T;E)$ fitting $m_V$ and including a forest $F$ with $d_F(t)=
m_{\rm for}(t)$ \ $(t\in T)$ if and only if \eqref{(ujterm.GR)} holds
and \begin{equation}\widetilde m_S(X) + \widetilde m_T(Y) -\vert X\vert \vert
Y\vert + \widetilde m_{\rm for} (T-Y) - \vert T-Y\vert -\vert X\vert
+1 \leq \gamma \ \hbox{ whenever }\ \emptyset \not =X\subseteq S, \
Y\subseteq T. \label{(T2fax)} \end{equation} \end{theorem}

\medskip We also remark that the results of Section \ref{poli} can be
used in a similar way to generalize Theorem \ref{wooded} so as to have
upper and lower bounds for the degrees of the nodes.

\subsubsection{Wooded hypergraphs}

A hypergraph is called {\bf wooded} if it can be trimmed to a graph
which is a forest, that is, if it is possible to select two distinct
elements from each hyperedge in such a way that the selected pairs, as
graph edges, form a forest.  Suppose we have a hypergraph $H=(S, {\cal
T})$ on node-set $S$.  It is well known that $H$ can be represented
with a simple bipartite graph $G_H=(S,T;E)$ where the elements of $T$
correspond to the hyperedges and the set of neighbours of $t\in T$ in
$G_H$ is just the hyperedge corresponding to $t$.  Obviously, $H$ is
wooded precisely if the associate bipartite graph $G_H$ has a
$T_2$-forest.  In this terminology, Theorem \ref{Lovasz1m} asserts
that a hypergraph is wooded if and only if the union of any $j>0$
hyperedges has at least $j+1$ elements.

Theorem \ref{wooded} can also be reformulated in terms of wooded
hypergraphs but here we do this only for the special case when the
hypergraph is $\ell$-uniform where $\ell\geq 2$.

\begin{corollary}\label{luniform} Let $m_S$ be a degree-specification on $S$
with $\widetilde m_S(S)=\gamma $ and let $\ell\geq 2$ be an integer.
There is an $\ell$-uniform wooded hypergraph fitting $m_S$ if and only
if $\tau := \gamma / \ell $ is an integer and \begin{equation}m_S(s) \leq \tau
\leq \vert S\sp +\vert -1 \ \hbox{for}\ s\in S\sp + \label{(elerdos)}
\end{equation} where $S\sp +=\{s\in S:  m_S(s)>0\}$.  \end{corollary}

\proof{Proof.} As the necessity of the conditions is straightforward, we
consider only sufficiency.  Since nodes $s\in S$ with $m_S(s)=0$ will
not belong to any hyperedge, we can delete them, and thus assume that
$S\sp +=S$.  Note that \eqref{(elerdos)} implies that $\widetilde
m_S(X)\leq \tau \vert X\vert $ for every $X\subseteq S$.

Let $T$ be a set of $\tau $ new elements.  Define $m_T(t):=\ell $ for
each $t\in T$ and let $p_T$ be a set-function on $T$ defined in
\eqref{(termrank.pTdef1)}.  If there is a simple bigraph $G=(S,T;E)$
covering $p_T$ and fitting $(m_S,m_T)$, then $G$ has a $T_2$-forest
and the hypergraph on $S$ associated with $G$ is an $\ell$-uniform
wooded hypergraph, in which case we are done.

Suppose that the requested bigraph does not exist.  Then one of the
conditions in Theorem \ref{wooded} fails to hold.  Suppose first that
there are sets $X\subseteq S, Y\subseteq T$ violating
\eqref{(ujterm.GR)}, that is, $ \widetilde m_S(X) + \widetilde m_T(Y)
-\vert X\vert \vert Y\vert > \gamma ,$ implying $$ \widetilde m_S(X) +
\vert Y\vert (\ell-\vert X\vert ) = \widetilde m_S(X) + \ell\vert
Y\vert -\vert X\vert \vert Y\vert > \gamma = \tau \ell .$$ If
$\ell\geq \vert X\vert $, then $$ \widetilde m_S(X) + \tau (\ell-\vert
X\vert )= \widetilde m_S(X)+ \vert T\vert (\ell-\vert X\vert ) \geq
\widetilde m_S(X)+ \vert Y\vert (\ell-\vert X\vert ) >\tau \ell,$$
from which $\widetilde m_S(X)> \tau \vert X\vert $, a contradiction.

If $\ell<\vert X\vert $, then $$ \gamma =\widetilde m_S(S) \geq
\widetilde m_S(X) \geq \widetilde m_S(X) + \vert Y\vert (\ell-\vert
X\vert ) >\gamma ,$$ a contradiction again, showing that
\eqref{(ujterm.GR)} holds.

Consider now the case when there are sets $\emptyset \not =X\subseteq
S, \ Y\subseteq T$ violating \eqref{(T2fa)}, that is, $$ \widetilde
m_S(X) + \widetilde m_T(Y) -\vert X\vert \vert Y\vert -\vert X\vert
-\vert Y\vert + \tau +1 >\gamma $$ from which $$ \widetilde m_S(X)+
\vert Y\vert (\ell-\vert X\vert -1) - \vert X\vert + \tau +1 =
\widetilde m_S(X) + \ell\vert Y\vert -\vert X\vert \vert Y\vert -\vert
X\vert -\vert Y\vert + \tau +1 >\gamma =\tau \ell.$$

If $\ell >\vert X\vert +1$, then $$ \widetilde m_S(X)+ \tau
(\ell-\vert X\vert -1) - \vert X\vert + \tau +1\geq \widetilde m_S(X)+
\vert Y\vert (\ell-\vert X\vert -1) - \vert X\vert + \tau +1 >\tau
\ell$$ from which $ \widetilde m_S(X) - \tau \vert X\vert -\vert
X\vert +1 >0,$ and hence $\tau \vert X\vert \geq \widetilde m_S(X) >
\tau \vert X\vert + \vert X\vert -1,$ implying that $X=\emptyset $, a
contradiction.

Suppose now that $\ell \leq \vert X\vert +1$.  Since $m_S(s)$ is
positive for every $s\in S$, we have $\widetilde m_S(S)-\vert S\vert
\geq \widetilde m_S(X)-\vert X\vert $. Hence $$ \widetilde m_S(S) -
\vert S\vert + \tau +1 \geq \widetilde m_S(X)-\vert X\vert + \tau +
1\geq \widetilde m_S(X)+ \vert Y\vert (\ell-\vert X\vert -1) - \vert
X\vert + \tau +1 >\tau \ell$$ from which $$ \tau \ell- \vert S\vert +
1 = \widetilde m_S(S) - \vert S\vert + 1 > \tau \ell - \tau ,$$ that
is, $\tau >\vert S\vert -1$, contradicting \eqref{(elerdos)}.  $\bullet
$

\medskip

T. Kir\'aly \cite{TKiraly} pointed out that there is a simple direct
proof of Corollary \ref{luniform} not relying on the theory of
supermodular functions.

In \cite{Berczi-Frank16b}, we describe two other extensions of Ryser's
theorem.

\medskip \medskip

\noindent {\bf Acknowledgements} \ \ We are grateful to Richard
Brualdi for the correspondence in which he always promptly provided us
with extremely useful background information on the history of the
topic.  We also thank Lap Chi Lau for his remarks.  Special thanks are
due to Zolt\'an Szigeti for carefully checking the details of a first
draft.  The two anonymous referees provided a great number of
particularly useful comments.  We gratefully acknowledge their
precious efforts.

The research was supported by the Hungarian Scientific Research Fund -
OTKA, No K109240.  The work of the first author was financed by a
postdoctoral fellowship provided by the Hungarian Academy of Sciences.

\medskip 

\end{document}